\documentclass[11pt]{article}
\usepackage{amssymb, amsmath,  amsfonts,dsfont, mathrsfs,euscript,eufrak,colonequals}
\usepackage{textcomp}
\usepackage{latexsym}
\usepackage{indentfirst}
\usepackage[normalem]{ulem}

\newtheorem{Theorem}{Theorem}

\newtheorem{Corollary}{Corollary}
\newtheorem{Remark}{Remark}

\newcommand{\defeq}{\colonequals}
\newcommand{\sgn}{\mathrm{sgn}}
\newcommand{\dd}{\,d}
\newcommand{\id}{\mathds{1}}

\newcommand{\R}{\mathbb R}

\newcommand{\N}{\mathbb N}

\newcommand{\Imagpart}{\mathrm{Im}\,}

\newcommand{\BV}{\boldsymbol{V}}

\newcommand{\ID}{\boldsymbol{I}}
\newcommand{\Q}{\boldsymbol{Q}}

\newcommand{\e}{\varepsilon}

\newcommand{\Ln}{\mathrm{Ln}}
\newcommand{\Arg}{\mathrm{Arg}}

\begin{document}
\author{A. A. Khartov$^{1,2,3,}$\footnote{Email addresses: \texttt{khartov.a@iitp.ru}, \texttt{alexeykhartov@gmail.com}}}
\title{Some criteria of rational-infinite divisibility for probability laws}

\footnotetext[1]{Institute for Information Transmission Problems (Kharkevich Institute) of Russian Academy of Sciences, Bolshoy Karetny per. 19, build.1, 127051 Moscow, Russia.}
\footnotetext[2]{Laboratory for Approximation Problems of Probability, Smolensk State University, 4 Przhevalsky st., 214000 Smolensk, Russia. }
\footnotetext[3]{Saint-Petersburg National Research University of Information Technologies, Mechanics and Optics (ITMO University), 49 Kronverksky Pr., 197101 Saint-Petersburg, Russia.}

\maketitle
\begin{abstract}
	We study the class $\Q$ of distribution functions $F$ that have the property of rational-infinite divisibility: there exist some infinitely divisible distribution functions $F_1$ and $F_2$ such that $F_1=F*F_2$. The class $\Q$  is a wide natural extension of the fundamental class of  infinitely divisible distribution functions. We are interested in general conditions  to belong to the class $\Q$ in terms of characteristic functions. We obtain criteria that  seem to be convenient for the application for some cases, and we illustrate it by several examples in the paper.
\end{abstract}


\section{Introduction}
Let $F$ be a  distribution function on the real line $\R$ with the characteristic function 
\begin{eqnarray*}
	f(t)\defeq\int_{\R} e^{itx} \dd F(x),\quad t\in\R.
\end{eqnarray*}
Recall that $F$ (and the corresponding law) is called \textit{infinitely divisible} if for every positive integer $n$ there exists a  distribution function $F_{1/n}$ such that $F=(F_{1/n})^{*n}$, where ``$*$'' denotes a convolution operation, i.e. $F$ is $n$-fold convolution power of $F_{1/n}$.  It is well known that $F$ is infinitely divisible if and only if the characteristic function $f$ is represented by \textit{the L\'evy--Khinchine formula}:
\begin{eqnarray}\label{repr_f}
	f(t)=\exp\biggl\{it \gamma+\int_{\R} \bigl(e^{itx} -1 -it \sin(x)\bigr)\tfrac{1+x^2 }{x^2} \dd G(x)\biggr\},\quad t\in\R,
\end{eqnarray}
with some \textit{shift parameter} $\gamma\in\R$, and with a  bounded non-decreasing \textit{spectral function} $G: \R \to \R$  that is assumed to be right-continuous at every point of the real line with the condition $G(-\infty)=0$. We use $x\mapsto \sin (x)$ as the ``centering function'' in the integral in \eqref{repr_f} following to Zolotarev \cite{Zolot1} and \cite{Zolot2}, but  there are also some other variants (see \cite{Sato1999}, p. 38).  It is well known that \textit{the spectral pair} $(\gamma, G)$ is uniquely determined by $f$ and hence by $F$. The L\'evy--Khinchine formula is  very important for probability theory and for it's  applications in related fields (see \cite{Appl} and \cite{Sato}).

There is a natural extension for the class of infinitely divisible distribution functions. We call a distribution function $F$ (and the corresponding law) \textit{rational-infinitely divisible} if there exist some infinitely divisible distribution functions $F_1$ and $F_2$ such that $F_1=F*F_2$. This equality may be written in terms of characteristic functions: $f(t)=f_1(t)/f_2(t)$, $t\in\R$, where $f_1$ and $f_2$ denote characteristic functions of $F_1$ and $F_2$, respectively. It is not difficult to show in this case that $f(t)\ne 0$ for any $t\in\R$, and $f$ admits representation \eqref{repr_f} with $\gamma=\gamma_1-\gamma_2$ and $G=G_1-G_2$, where $(\gamma_1, G_1)$ and $(\gamma_2, G_2)$ denote the spectral pairs of $F_1$ and $F_2$, respectively. Note that the function $G$ is of bounded variation on $\R$ (it is non-monototic in general), right-continuous at every point and $G(-\infty)=0$. The class of all functions satisfying these properties will be denoted by $\BV$. The  pair $(\gamma, G)$ is uniquely determined by $f$ and by $F$  as for infinitely divisible distribution functions (see \cite{Khinch} p. 31-32, or \cite{GnedKolm} p. 80). Now suppose that $f$ is represented by formula \eqref{repr_f} with some $\gamma\in\R$ and  $G\in\BV$. Following Lindner and Sato \cite{LindSato}, the corresponding distribution function $F$  is called \textit{quasi-infinitely divisible}. Due to the Hahn--Jordan decomposition for $G$, it is not difficult to show that $F$ is rational-infinitely divisible. The  first detailed analysis of the class of such distribution functions on $\R$ was performed in \cite{LindPanSato}. There are interesting applications in theory of stochastic processes (see  \cite{LindSato} and \cite{Pass}), number theory (see \cite{Nakamura}),  physics (see \cite{ChhDemniMou} and \cite{Demni}), and insurance mathematics (see  \cite{ZhangLiuLi}). 
  
The class of all rational-infinitely divisible distribution functions we will denote by $\Q$. By definition, it is clear  that this class contains the class $\ID$ of all infinitely divisible distribution functions and, moreover,  it is essentially wider  than the latter. So $\Q$ includes all discrete distribution functions $F$, whose characteristic functions $f$ are separeted from zero, i.e.  $|f(t)|\geqslant \mu$  for some constant $\mu>0$ and for any $t\in\R$ (see \cite{AlexeevKhartov} and \cite{Khartov2}). The class $\Q$ is also contains mixtures of discrete and absolutely continuous distribution functions, where the discrete part is non-zero  and its characteristic function is separeted from zero (see \cite{Berger} and  \cite{BergerKutlu}).  Note that a rational-infinitely divisible distribution functions may  have a bounded set of growth points (support of distribution)  that is impossible for non-degenerate infinitely divisible distribution functions. There are also some interesting particular sufficient conditions  to belong to $\Q$ (see \cite{LindPanSato}). However, we think that the problem about general criteria of rational-infinite divisibility for probability laws is not completely solved. Indeed, for instance, there are no criteria that can be conveniently applied to  absolutely continuous laws. In particular, we are not aware any criterion (excluding definitions) that would take a characteristic function of a Gaussian law or a stable law (or another important law) as input and confirm on the output that it is corresponded to (rational-)infinitely divisible distribution function. In this work, we propose some results for  this problem. So the obtained criteria yield general methods to determine whether a given characteristic function corresponds to some rational-infinitely or purely infinitely divisible law. These criteria seem to be rather convenient and easy in applications for some cases, and we illustrate it  in this paper by special examples, where we test the criteria with characteristic functions of well known absolutely continuous laws. Of course, we don't asert that the proposed criteria are the best ways to solve the problem. Also it should be noted that we tried to formulate the criteria with priority on generality and hence, dealing with particular examples, their own additional analysis is usually needed. Some techniques for this are demonstrated in the mentioned examples.
 
Note that the proposal of using of the term ``rational-infinitely divisible laws'' for the described distributions is an attempt to reflect  the key decomposition property (that $F_1=F*F_2$ with some $F_1,F_2\in\ID$) in the class name similarly to that it was done for the important classes of stable, self-decomposable  or infinitely divisible laws. We think that a classification of probability laws based on the existence and various forms of the Levy--Khinchine representation may be not convenient in the perspective. So, in the paper \cite{KhartovLogarithm}, it was showed that an \textit{arbitrary} characteristic function $f$ admits a special L\'evy--Khinchine type representation on segments, where $f(t)\ne 0$. Also there exists probability laws whose characteristic function admit representation \eqref{repr_f} with  spectral function $G$ that has an \textit{infinite} total variation on $\R$. It is rather interesting direction,  but we will not discuss these questions here.

The statements of the results are given in Section~2, their proofs are  provided in Section~4, the examples of application are written in Section~3. Throughout the paper, we use the following notation. We denote by $\N$  the sets of positive integers and let $\N_0\defeq \N\cup \{0\}$. For every $G\in\BV$ its total variation on $\R$ will be denoted by $\|G\|$ and the total variation on $(-\infty,x]$ --- by $|G|(x)$, $x\in\R$. So we have $|G(x)|\leqslant |G|(x)\leqslant\|G\|$, $x\in\R$, and $|G|(+\infty)=\|G\|$.   Let $\lfloor \cdot \rfloor$ be the floor function, i.e. for any $x\in\R$ we have $\lfloor x \rfloor= k$  with integer $k$ such that $k\leqslant x<k+1$. Below $\id_a$ with fixed $a\in\R$ denotes the following function: $\id_a(x)=1$ for $x\geqslant a$ and $\id_a(x)=0$ for $x<a$. We always set for the sums $\sum_{j=n}^m a_j=0$ if $n>m$.

\section{General results}

We now only assume that $f$ is a complex-valued  continuous function on $\R$ satisfying $f(0)=1$ and $f(t)\ne 0$ for any $t\in\R$. We define \textit{the distinguished logarithm} $\Ln f$, i.e. $\Ln f(t):=\ln |f(t)|+i \Arg f(t)$, $t\in\R$, where the argument $\Arg f$ is uniquely defined on $\R$ by continuity with the condition $\Arg f(0)=0$. We next introduce the following function
\begin{eqnarray}\label{def_psi}
	\psi(t,h):= \Ln f(t)-\dfrac{1}{2}\bigl(\Ln f(t-h)+\Ln f(t+h)\bigr),\quad t, h\in\R.
\end{eqnarray}
Observe that 
\begin{eqnarray}\label{def_Delta2Ln}
	\Delta^2_h\Ln f(t):=-2 \psi(t,h)=\Ln f(t-h)+\Ln f(t+h)-2\Ln f(t),\quad t, h\in\R,
\end{eqnarray}
is a finite difference of the second order for $\Ln f$. It should be noted that the function  $\psi$ is often rather simpler and easier handling than $\Ln f$. It is also used for establishing some properties of $f$ (see \cite{Khartov2} p. 3, \cite{Khartov3} p. 119, and \cite{Petr} p. 34). The following theorem  shows that $\psi$ contains in fact full information about $f$.

\begin{Theorem}\label{th_crit_psi}
Suppose that 
\begin{eqnarray}\label{def_Idelta}
	\psi(t,h)= I(t,h)+\delta(t,h)\quad\text{for any}\quad t\geqslant0,\,\,h>0,
\end{eqnarray}
where the function $I$ admits the following representation
\begin{eqnarray}\label{def_I}
	I(t,h)= \int_{\R}e^{itx}\bigl(1-\cos(h x)\bigr)\tfrac{1+x^2}{x^2}\dd G(x)\quad\text{for any}\quad t\geqslant0,\,\, h>0,
\end{eqnarray}
with some $G\in\BV$, and the function $\delta$ satisfies the following condition:
	\begin{eqnarray}\label{th_crit_psi_cond_a}
		\sum_{k=0}^{n_{t,l}-1}\bigl(n_{t,l}-k\bigr)|\delta(kh_l,h_l)|\to 0, \quad l\to\infty,\quad\text{for any}\quad t>0,
	\end{eqnarray}
where  $(h_l)_{l\in\N}$ is a fixed sequence of positive reals such that $h_l\to 0$, $l\to\infty$, and $n_{t,l}:= \lfloor t/ h_l\rfloor$, $t>0$, $l\in\N$. Then  $f$ admits representation \eqref{repr_f} with the function $G$ and  $\gamma=\Arg f(1)$. If it is additionally known that $f$ is the characteristic function of some distribution function $F$, then $F\in\Q$.
\end{Theorem}

Here and below we use the following convention
\begin{eqnarray}\label{eq_conven_h2}
	\Bigl[\bigl( 1- \cos(h x)\bigr)\tfrac{1+x^2}{x^2}\Bigr]\biggl|_{x=0}\defeq\lim\limits_{x\to 0}\Bigl[\bigl( 1- \cos(h x)\bigr)\tfrac{1+x^2}{x^2}\Bigr]=\dfrac{h^2}{2}.
\end{eqnarray}

\begin{Corollary}\label{co_th_crit_psi}
	Suppose that for any $t\geqslant0$ and  $h>0$ the formula
	\begin{eqnarray}\label{co_th_crit_psi_eq}
		\psi(t,h)= \int_{\R}e^{itx}\bigl(1-\cos(hx)\bigr)\tfrac{1+x^2}{x^2}\dd G(x)
	\end{eqnarray}
	holds with some $G\in\BV$. Then \eqref{th_crit_psi_cond_a} obviously holds $($with $\delta(kh,h):=0$ for all $k\in\N_0$ and $h>0${}$)$,  and $f$ admits representation \eqref{repr_f} with the function $G$ and $\gamma=\Arg f(1)$. If it is additionally known that $f$ is the characteristic function of some distribution function $F$, then $F\in\Q$.
\end{Corollary}

We note that the representation by formulas \eqref{def_Idelta} or \eqref{co_th_crit_psi_eq} is not only sufficient for formula \eqref{repr_f}, but it is also  necessary.

\begin{Remark}
	Suppose that $f$ admits representation \eqref{repr_f}  with some $\gamma\in\R$ and $G\in\BV$ $($\,for instance,  $f$ is the characteristic function of some distribution function $F\in\Q$\,$)$. Then
	\begin{eqnarray}\label{th_crit_psi_eq}
		\psi(t,h)=\int_{\R}e^{itx}\bigl(1-\cos(hx)\bigr)\tfrac{1+x^2}{x^2}\dd G(x)\quad\text{for all}\quad t,h\in\R.
	\end{eqnarray}
\end{Remark}
Indeed, substituting \eqref{repr_f} into \eqref{def_psi}, for any $t,h\in\R$ we get
\begin{eqnarray*}
	\psi(t,h)&=&it \gamma+\int_{\R} \Bigl(e^{itx} -1 -it \sin(x)\Bigr)\tfrac{1+x^2 }{x^2}  \dd G(x)\nonumber\\
	&&{}-  i\,\tfrac{t-h+t+h}{2}\,\gamma-\int_{\R} \Bigl(e^{itx}\,\tfrac{e^{-ihx}+e^{ihx}}{2}  - i\,\tfrac{t-h+t+h}{2}\,\sin(x)\Bigr)\tfrac{1+x^2 }{x^2}  \dd G(x)\nonumber\\
	&=&\int_{\R} e^{itx} \bigl(1-\cos(hx)\bigr) \tfrac{1+x^2}{x^2}\dd G(x).
\end{eqnarray*}

Sometimes it is  more convenient to work with the following functions
\begin{eqnarray}\label{def_phi}
	\varphi_{\pm}(t,h):=e^{\pm 2\psi(t,h)}-1=\biggl(\dfrac{f(t)^2}{f(t-h)f(t+h)}\biggr)^{\pm 1}-1,\quad t, h\in\R.
\end{eqnarray}
Using these functions, one can formulate the sufficient condition for representation \eqref{repr_f}.

\begin{Theorem}\label{th_suff_phi}
	 Suppose that 
	 \begin{eqnarray}\label{th_suff_phi_bkh}
	 	\varphi_{\pm}(t,h)=\pm 2I(t,h)+\rho(t,h)\quad\text{for any}\quad t\geqslant0,\,\,h>0,
 	 \end{eqnarray}
	 where the function $I$ admits representation \eqref{def_I} with some $G\in\BV$, the function $\rho$ satisfies the following condition:
	\begin{eqnarray}\label{th_suff_phi_b}
	\sum_{k=0}^{n_{t,l}-1}\bigl(n_{t,l}-k\bigr)|\rho(kh_l,h_l)|\to 0, \quad l\to\infty,\quad\text{for any}\quad t>0.
	\end{eqnarray}
Also suppose that
\begin{eqnarray}\label{th_suff_phi_phi2}
	\sum_{k=0}^{n_{t,l}-1}\bigl(n_{t,l}-k\bigr)\bigl|\varphi_{\pm}(kh_l,h_l)\bigr|^2\to 0, \quad l\to\infty,\quad\text{for any}\quad t>0,
\end{eqnarray}
where $(h_l)_{l\in\N}$ is a fixed sequence of positive reals such that $h_l\to 0$, $l\to\infty$, and $n_{t,l}:= \lfloor t/ h_l\rfloor$, $t>0$.	Then  $f$ admits representation \eqref{repr_f} with the function $G$ and  $\gamma=\Arg f(1)$. If it is additionally known that $f$ is the characteristic function of some distribution function $F$, then $F\in\Q$.
\end{Theorem}

Let us return to the function \eqref{def_psi}. We saw that the admission of representation  \eqref{th_crit_psi_eq} is a criterion for rational-infinite divisibility of distribution functions. Let us write it in the following form:
\begin{eqnarray*}
	\dfrac{\Delta^2_h\Ln f(t)}{h^2}=-\int_{\R}e^{itx}\cdot\dfrac{\sin^2\bigl(\tfrac{hx}{2}\bigr)}{\bigl(\tfrac{hx}{2}\bigr)^2}\cdot(1+x^2)\dd G(x).
\end{eqnarray*}
If $h\to 0$ then the left-hand side  tends to $(\Ln f)''(t)$, when the latter exists, and the fraction  formally dissapears  from the integrand in the right-hand side. This idea leads us  to the following result.

\begin{Theorem}\label{th_SecDerivLnf}
	$a)$ Suppose that $f$ admits representation \eqref{repr_f} with some $\gamma\in\R$  and $G\in\BV$ $($\,for instance,  $f$ is the characteristic function of some distribution function $F\in\Q$\,$)$ satisfying
	\begin{eqnarray}\label{th_SecDerivLnf_int}
		\int_{\R} (1+x^2)\dd |G|(x)<\infty.
	\end{eqnarray}
	Then for any $t\in \R$ there exists $(\Ln f)''(t)$ that admits the following representation	
	\begin{eqnarray}\label{th_SecDerivLnf_2deriv}
		(\Ln f)''(t)= -\int_{\R} e^{itx}(1+x^2)\dd G(x),\quad t\in\R.
	\end{eqnarray}
	$b)$ Suppose that the formula \eqref{th_SecDerivLnf_2deriv}  holds with some function $G\in\BV$ satisfying \eqref{th_SecDerivLnf_int}. Then  $f$ admits representation \eqref{repr_f} with the function $G$ and  $\gamma=\Arg f(1)$. If it is additionally known that $f$ is the characteristic function of some distribution function $F$, then $F\in\Q$.
\end{Theorem}

Theorems \ref{th_crit_psi}--\ref{th_SecDerivLnf}  can be applied as methods of determining that a given function $f$ is the characteristic for some infinitely divisible law. The following assertion, which we formulate as theorem, is a simple corollary from the previous theorems and the mentioned facts  about infinitely divisible laws (see Introduction).

\begin{Theorem}\label{th_InfDiv}
		Suppose that the assumptions of Theorem \ref{th_crit_psi}, or Corollary \ref{co_th_crit_psi}, or Theorem \ref{th_suff_phi}, or Theorem \ref{th_SecDerivLnf} $b)$ hold with some non-decreasing function $G\in\BV$. Then $f$ is a characteristic function of some infinitely divisible distribution function $F$, i.e. $F \in \ID$.
\end{Theorem}

\section{Examples of application}
In this section, we illustrate application  of Theorem \ref{th_crit_psi}--\ref{th_InfDiv} for characteristic functions of known infinitely divisible probability laws. Note that, of course, we will not use this knowledge in analysis to obtain  the conclusions by the criteria. We always start the analysis with a given function $f$ without any additional information. We deal with known distributions in order to show clearly  the correctness of the results of application.\\

\textbf{Example 1.} Here we show that the general results from the previous section can be formally applied to the characteristic function of the Gaussian law:
\begin{eqnarray*}
	f(t)=\exp\biggl\{it\gamma-\dfrac{\sigma^2 t^2}{2} \biggr\},\quad t\in\R,
\end{eqnarray*}
where $\gamma\in\R$ and $\sigma^2\geqslant 0$ (we include the degenerate case $\sigma^2=0$).

Since $f(t)\ne 0$, $t\in\R$, and $f(0)=1$, we have
\begin{eqnarray*}
	\Ln f(t)=it\gamma-\dfrac{\sigma^2 t^2}{2},\quad t\in\R.
\end{eqnarray*}
The easiest way is to use Theorem \ref{th_SecDerivLnf}. So  $(\Ln f)''(t)=-\sigma^2$ for any $t\in\R$. It is easily seen that $(\Ln f)''$  admits representation \eqref{th_SecDerivLnf_2deriv} if we set $G(x):= \id_{\sigma^2}(x)$, $x\in\R$ (in particular, $G(x)\equiv0$   for the case $\sigma^2=0$).  The function $G\in\BV$ satisfies \eqref{th_SecDerivLnf_int} and it is non-decreasing on $\R$, $\gamma=\Arg f(1)$. By Theorem~\ref{th_SecDerivLnf} and Theorem~\ref{th_InfDiv}, we conclude that $f$ is the characteristic function of an infinitely divisible distribution function $F$ with the spectral pair $(\gamma,G)$.

Also it is not difficult to apply another criterion. Let us consider the function $\psi$:
\begin{eqnarray*}
	\psi(t,h)
	&=&\Ln f(t)-\dfrac{1}{2}\bigl(\Ln f(t-h)+\Ln f(t+h)\bigr)\\
	&=&it\alpha-\dfrac{\sigma^2 t^2}{2}- \dfrac{1}{2}\biggl( i(t-h)\alpha-\dfrac{\sigma^2 (t-h)^2}{2} + i(t+h)\alpha-\dfrac{\sigma^2 (t+h)^2}{2} \biggr)\\
	&=&it\alpha-\dfrac{\sigma^2 t^2}{2}- \dfrac{1}{2}\bigl( 2it\alpha-\sigma^2t^2-\sigma^2 h^2  \bigr)= \dfrac{\sigma^2 h^2}{2},\quad t,h\in\R.
\end{eqnarray*}
On account of the convention in \eqref{eq_conven_h2}, we can write
\begin{eqnarray*}
	\psi(t,h)=\dfrac{\sigma^2 h^2}{2}=\int_{\R}e^{itx}\bigl(1-\cos(hx)\bigr)\tfrac{1+x^2}{x^2}\dd G(x),\quad t,h\in\R,
\end{eqnarray*}
with $G(x)=\id_{\sigma^2}(x)$, $x\in\R$. According to Theorem \ref{th_crit_psi} with Corollary \ref{co_th_crit_psi} and Theorem~\ref{th_InfDiv}, we come to the same conclusion.\quad $\Box$\\

\textbf{Example 2.}  Let us consider the function
\begin{eqnarray*}
	f(t)=\dfrac{1}{\cosh (t)},\quad t\in\R.
\end{eqnarray*}
We first show that $f$ is the characteristic function of absolutely continuous probability law.  Since $f\in L_1(\R)$, we may define the function
\begin{eqnarray*}
	p(x):= \dfrac{1}{2\pi} \int_{\R} e^{-itx} f(t)\dd t= \dfrac{1}{2\pi} \int_{\R}\dfrac{ e^{-itx}}{\cosh (t)} \dd t=\dfrac{1}{\pi} \int_{0}^\infty\dfrac{\cos(tx)}{\cosh(t)} \dd t,\quad x\in\R.
\end{eqnarray*}
The last equality is valid, because $\cosh (\cdot)$ is an even  function. Next, from \cite{GradRyz}, formula \textbf{3.981} 3., we know that
\begin{eqnarray*}
	\int_{0}^\infty\dfrac{\cos(as)}{\cosh(\beta s)} \dd s= \dfrac{\pi}{2\beta\cosh \Bigl(\dfrac{\pi a}{2\beta}\Bigr)}\quad\text{for any}\quad a\in\R,\,\, \beta>0.
\end{eqnarray*}
Therefore
\begin{eqnarray*}
	p(x)=  \dfrac{1}{2\cosh \Bigl(\dfrac{\pi x}{2}\Bigr)},\quad x\in\R.
\end{eqnarray*}
It is clear that $p(x)\geqslant 0$ for any $x\in\R$ and $p\in L_1(\R)$. Then we may write
\begin{eqnarray*}
	f(t)=\int_{\R} e^{itx} p(x) \dd x,\quad t\in\R,
\end{eqnarray*}
and, in particular, 
\begin{eqnarray*}
	\int_{\R} p(x) \dd x=f(0)=\dfrac{1}{\cosh(0)}=1.
\end{eqnarray*}
Thus $p$ is the  density of absolutely continuous distribution with function $F(x)=\int_{-\infty}^{x} p(u)\dd u$, $x\in\R$, and $f$ is its characteristic function.

We now show by Theorem \ref{th_SecDerivLnf} and Theorem \ref{th_InfDiv} that $F$ is infinitely divisible. We write for any $t\in\R$:
\begin{eqnarray*}
	\Ln f(t)=-\ln \cosh(t),\quad (\Ln f)'(t)=-\tanh(t),\quad (\Ln f)''(t)=-\dfrac{1}{\cosh^2(t)}.
\end{eqnarray*}
Thus $(\Ln f)''(t)=-f(t)^2$, $t\in\R$. The function $t\mapsto f(t)^2$ is characteristic for the distribution with the density
\begin{eqnarray*}
	p^{*2}(x)=\int_{\R}p(x-y)p(y)\dd y,\quad x\in\R.
\end{eqnarray*}
Hence
\begin{eqnarray*}
	(\Ln f)''(t)=-\int_{\R} e^{itx} p^{*2}(x) \dd x=-\int_{\R} e^{itx} (1+x^2)\dd G(x),\quad t\in\R,
\end{eqnarray*}
where we set
\begin{eqnarray*}
	G(x):=\int_{-\infty}^x \dfrac{p^{*2}(u)\dd u}{1+u^2},\quad x\in\R.
\end{eqnarray*}
Due to the summability and the non-negativity of the integrand function here, $G\in \BV$ and it is non-decreasing on $\R$. The condition \eqref{th_SecDerivLnf_int} is satisfied.  Thus, by Theorems \ref{th_SecDerivLnf} and \ref{th_InfDiv}, $F$ is infinitely divisible distribution function with the spectral components $\gamma=\Arg f(1)=0$ and  $G$.

Actually, it has been already known that $f$ corresponds to an infinitely divisible distribution. In \cite{Lukacs} (see p. 88), this fact was proved by more difficult way, namely, using the Weierstrass factorization theorem and the continuity theorem (for special compositions of the Laplace distributions). \quad $\Box$\\

\textbf{Example 3.}  Here we deal with the characteristic function of the general Cauchy distribution with asymmetry (general continuous stable law with $\alpha=1$):
\begin{eqnarray}\label{def_f_Cauchy}
	f(t)=\exp\biggl\{it\gamma-\lambda|t|-i\dfrac{2}{\pi}\,\lambda\beta\, t \ln |t|\biggr\},\quad t\in\R,
\end{eqnarray}
where $\gamma\in\R$, $\lambda>0$, and $\beta\in[-1,1]$. 

So $f(t)\ne 0$ for any $t\in\R$, $f(0)=1$, and
\begin{eqnarray*}
	\Ln f(t)=it\gamma-\lambda|t|-i\dfrac{2}{\pi}\,\lambda\beta\, t \ln |t|,\quad t\in\R.
\end{eqnarray*}
Let us write the function $\psi$:
\begin{eqnarray}
	\psi(t,h)&:=&\Ln f(t)-\dfrac{1}{2}\bigl(\Ln f(t-h)+\Ln f(t+h)\bigr)\nonumber\\
	&=& it\gamma-\lambda|t|-i\dfrac{2}{\pi}\,\lambda\beta\, t \ln |t|\nonumber\\
	&&{}-\dfrac{1}{2}\biggl( i(t-h)\gamma-\lambda|t-h|-i\dfrac{2}{\pi}\,\lambda\beta\, (t-h)\ln |t-h|\nonumber\\
	&&{}\,\,\quad+i(t+h)\gamma-\lambda|t+h|-i\dfrac{2}{\pi}\,\lambda\beta\, (t+h)\ln |t+h|\biggr)\nonumber\\
	&=&\lambda\Delta_{h}(t)+i\,\dfrac{\lambda\beta}{\pi}\, A_h(t),\quad t\in\R,\,\, h\in\R,\label{eq_psiDeltahAh}
 \end{eqnarray}
where for any $t\in\R$ and $h\in\R$ we denote
\begin{eqnarray*}
	\Delta_h(t)&:=&\dfrac{1}{2}\,\bigl(|t-h|+|t+h|-2|t| \bigr),\\
	A_h(t)&:=&(t-h)\ln|t-h|+(t+h)\ln|t+h|-2t\ln|t|.
\end{eqnarray*}
Let us fix any $h>0$. We will show applying Fourier transforms that the functions $\Delta_h$ and $A_h$ admit suitable integral representations.

We first consider the triangular function $\Delta_h$:
\begin{eqnarray*}
	\Delta_h(t)=\begin{cases}
		t+h,& t\in(-h,0],\\
		h-t,& t\in(0,h),\\
		0,& |t|\geqslant h.
	\end{cases}
\end{eqnarray*}
We introduce
\begin{eqnarray}\label{def_Dh}
	D_h(x):= \dfrac{1}{2\pi} \int_{\R} e^{-itx}  \Delta_h(t)\dd t= \dfrac{1}{\pi} \int_{0}^h \cos(tx) (h-t)\dd t= \dfrac{1-\cos(hx)}{\pi x^2},\quad x\in\R.
\end{eqnarray}
Here the last expression was directly calculated by the integration by parts. Since the functions $\Delta_h$ and $D_h$ are from $L_1(\R)$, we have
\begin{eqnarray}\label{eq_Deltah}
	\Delta_h(t)=\int_{\R} e^{itx} D_h(x) \dd x= \int_{\R} e^{itx} \dfrac{1-\cos(hx)}{\pi x^2} \dd x,\quad t\in\R.
\end{eqnarray}

We now turn to the function $A_h$. Observe that it is an odd function, and we consider it only for $t\geqslant0$:
\begin{eqnarray*}
	A_h(t)=(t-h)\ln|t-h|+(t+h)\ln(t+h)-2t\ln t,
\end{eqnarray*}
where we set $y \ln y:=0$ for $y=0$. The function $A_h$ goes to zero as $t\to\infty$. Indeed,
\begin{eqnarray*}
	A_h(t)&=& t\bigl(\ln(t-h)+\ln(t+h)-2\ln t \bigr) +h \bigl(\ln(t+h)-\ln(t-h)\bigr)\\
	&=& t \ln \bigl(1-(h/t)^2\bigr) + h\bigl(\ln \bigl(1+h/t\bigr)-\ln \bigl(1-h/t\bigr)\bigr),
\end{eqnarray*}
where, on account of the asymptotics  $\ln(1+y)\sim y$ as $y\to 0$, 
\begin{eqnarray*}
	\ln \bigl(1-(h/t)^2\bigr)\sim -(h/t)^2,\qquad \ln \bigl(1+h/t\bigr)-\ln \bigl(1-h/t\bigr)\sim 2h/t,\quad t\to\infty.
\end{eqnarray*}
Thus $A_h(t)\sim h^2/t\to 0$ as $t\to\infty$. Next, observe that $A_h$ is continuous on $[0,\infty)$  and, moreover, it has the derivatives at every $t>0$ and $t\ne h$:
\begin{eqnarray}\label{eq_AhDeriv}
	A_h'(t)=\ln|t-h|+\ln(t+h)-2\ln t, \quad  A_h''(t)=\dfrac{1}{t-h}+\dfrac{1}{t+h}-\dfrac{2}{t}.
\end{eqnarray} 
The function $A_h$ is decreasing on $(h,+\infty)$, because $A_h'(t)=\ln \bigl(1-(h/t)^2\bigr)<0$, $t>h$. 

On account of above-mentioned properties of $A_h$, we apply the Pringsheim  theorem about Fourier transforms for monotonic functions (see \cite{Fikhten3} p. 592, and \cite{Titch} p. 17) and we represent
\begin{eqnarray}\label{eq_AhFourier}
	A_h(t)=\lim_{R\to\infty}\dfrac{2}{\pi} \int_0^R \sin(tx) \biggl(\lim_{r\to\infty}\int_0^r \sin(sx) A_h(s)\dd s\biggr)\dd x,\quad t\geqslant 0.
\end{eqnarray}
Note that the formula is valid for $t=0$, because $A_h(0)=0$ (the Pringsheim theorem covers only the case $t>0$). Let us calculate the inner integral for $x>0$
\begin{eqnarray*}
	J_h(x)&:=& \lim_{r\to\infty}\int_0^r \sin(sx) A_h(s)\dd s\\
&=& \lim_{r\to\infty}\lim_{\e\to 0+}\Biggl[\int_0^{h-\e} \sin(sx) A_h(s)\dd s +\int_{h+\e}^r \sin(sx) A_h(s)\dd s\Biggr].
\end{eqnarray*}
We integrate by parts:
\begin{eqnarray*}
	J_h(x)&=& \lim_{r\to\infty}\lim_{\e\to 0+}\Biggl[  -\dfrac{\cos((h-\e)x)}{x}\,A_h(h-\e)+\dfrac{1}{x}\,A_h(0)  + \int_0^{h-\e} \dfrac{\cos(sx)}{x}\, A_h'(s)\dd s\\
	&&{}\qquad\qquad\,\,\, -\dfrac{\cos(rx)}{x}\,A_h(r)+\dfrac{\cos((h+\e)x)}{x}\,A_h(h+\e)+\int_{h+\e}^r \dfrac{\cos(sx)}{x}\, A_h'(s)\dd s\Biggr].
\end{eqnarray*}
Since the function $s\mapsto \tfrac{\cos(sx)}{x}\, A_h(s)$ is continuous, $A_h(0)=0$, and $A_h(r)\to 0$ as $r\to\infty$,  we get
\begin{eqnarray*}
	J_h(x)= \lim_{r\to\infty}\lim_{\e\to 0+}\Biggl[  \int_0^{h-\e} \dfrac{\cos(sx)}{x}\, A_h'(s)\dd s+\int_{h+\e}^r \dfrac{\cos(sx)}{x}\, A_h'(s)\dd s\Biggr].
\end{eqnarray*}
We again integrate by parts:
\begin{eqnarray*}
	J_h(x)&=& \lim_{r\to\infty}\lim_{\e\to 0+}\Biggl[ \dfrac{\sin((h-\e)x)}{x^2}\, A_h'(h-\e)-\lim_{s\to 0+}\biggl(\dfrac{\sin(sx)}{x^2}\, A_h'(s)\biggl) -\int_0^{h-\e} \dfrac{\sin(sx)}{x^2}\, A_h''(s)\dd s\\
	&&{}\qquad\qquad\,\,\,+ \dfrac{\sin(rx)}{x^2}\, A_h'(r)-\dfrac{\sin((h+\e)x)}{x^2}\, A_h'(h+\e)-\int_{h+\e}^r \dfrac{\sin(sx)}{x^2}\, A_h''(s)\dd s\Biggr].
\end{eqnarray*}
Here we have (see \eqref{eq_AhDeriv})
\begin{eqnarray*}
	\lim_{s\to 0+}\biggl(\dfrac{\sin(sx)}{x^2}\, A_h'(s)\biggl)&=&\lim_{s\to 0+}\biggl(\dfrac{\sin(sx)}{x^2}\, \bigl(\ln(h-s)+\ln(s+h)-2\ln s\bigr)\biggl)=0,\\
	\lim_{r\to \infty}\biggl( \dfrac{\sin(rx)}{x^2}\, A_h'(r)\biggr)&=&\lim_{r\to \infty}\biggl( \dfrac{\sin(rx)}{x^2}\,\ln \bigl(1-(h/r)^2\bigr)\biggr)= -\lim_{r\to \infty}\biggl( \dfrac{\sin(rx)h^2}{x^2r^2}\biggr)= 0.
\end{eqnarray*}
We next observe that
\begin{eqnarray*}
	A'(h\pm \e) &=&\ln|h\pm \e-h|+\ln(h\pm \e+h)-2\ln (h\pm \e)\\
	&=&\ln \e+\ln(2h\pm \e)-2\ln(h\pm \e)=\ln\e +\ln(2/h) +O(\e),\quad \e\to 0+.
\end{eqnarray*}
and $\sin((h\pm\e)x)=\sin (hx)+O(\e)$ as $\e\to 0+$. Hence
\begin{eqnarray*}
	\sin((h\pm\e)x)A'(h\pm \e)=\ln\e \cdot \sin(hx)+\ln(2/h)\cdot \sin(hx)+O(\e\ln \e),\quad \e\to 0+.
\end{eqnarray*}
Therefore
\begin{eqnarray*}
	\lim_{\e\to 0+}\biggl[ \dfrac{\sin((h-\e)x)}{x^2}\, A_h'(h-\e)-\dfrac{\sin((h+\e)x)}{x^2}\, A_h'(h+\e)\biggr]= \lim_{\e\to 0+} O(\e\ln \e)=0.
\end{eqnarray*}
Using \eqref{eq_AhDeriv}, we  next write :
\begin{eqnarray*}
	J_h(x)&=& \lim_{r\to\infty}\lim_{\e\to 0+}\Biggl[-\int_0^{h-\e} \dfrac{\sin(sx)}{x^2}\, A_h''(s)\dd s-\int_{h+\e}^r \dfrac{\sin(sx)}{x^2}\, A_h''(s)\dd s\Biggr]\\
	&=& -\dfrac{1}{x^2}\lim_{r\to\infty}\lim_{\e\to 0+}\Biggl[\int_0^{h-\e} \sin(sx) \biggl(\dfrac{1}{s-h}+\dfrac{1}{s+h}-\dfrac{2}{s}\biggr)\dd s\\
	&&{}\qquad\qquad\qquad\,\,\,+\int_{h+\e}^r \sin(sx) \biggl(\dfrac{1}{s-h}+\dfrac{1}{s+h}-\dfrac{2}{s}\biggr)\dd s\Biggr].
\end{eqnarray*}
We write
\begin{eqnarray*}
	J_h(x)&=&\dfrac{2}{x^2}\lim_{r\to\infty} \int_{0}^r\dfrac{\sin(sx)}{s}\,\dd s -\dfrac{2}{x^2}\lim_{r\to\infty}\lim_{\e\to 0+}\Biggl[\int_0^{h-\e}  \dfrac{s\sin(sx)}{s^2-h^2}\dd s+\int_{h+\e}^r\dfrac{s\sin(sx)}{s^2-h^2}\dd s\Biggr],
\end{eqnarray*}
where, since $x>0$, the first limit equals $\pi/2$ (see \cite{GradRyz} \textbf{3.721} 1.)  and the iterated limit equals $(\pi/2)\cdot \cos(hx)$ (see \cite{Fikhten2} p. 684, and  \cite{GradRyz} \textbf{3.723} 10.). Thus we obtain
\begin{eqnarray*}
	J_h(x)=\pi\, \dfrac{1-\cos(hx)}{x^2},\quad x>0.
\end{eqnarray*}
Let us return to formula \eqref{eq_AhFourier} and write:
\begin{eqnarray*}
	A_h(t)=\lim_{R\to\infty} 2 \int_0^R \sin(tx)\,  \dfrac{1-\cos(hx)}{x^2}\dd x,\quad t\geqslant 0.
\end{eqnarray*}
Since $x\mapsto \tfrac{1-\cos(hx)}{x^2}$ is an even function from $L_1(\R)$, we may rewrite this formula as follows:
\begin{eqnarray}\label{eq_AhInt}
	A_h(t)= \int_{\R} \sin(tx)\, \sgn(x)\, \dfrac{1-\cos(hx)}{x^2}\dd x=\dfrac{1}{i}\int_{\R} e^{itx}\, \sgn(x)\, \dfrac{1-\cos(hx)}{x^2}\dd x,\quad t\geqslant 0.
\end{eqnarray}

We now combine representations \eqref{eq_Deltah} and \eqref{eq_AhInt} in formula \eqref{eq_psiDeltahAh}:
\begin{eqnarray*}
	\psi(t,h)&=&\lambda\Delta_{h}(t)+i\,\dfrac{\lambda\beta}{\pi}\, A_h(t)\\
	&=&\lambda \int_{\R} e^{itx} \dfrac{1-\cos(hx)}{\pi x^2} \dd x  + i\,\dfrac{\lambda\beta}{\pi}\cdot \dfrac{1}{i} \int_{\R} e^{itx}\, \sgn(x)\, \dfrac{1-\cos(hx)}{x^2}\dd x\\
	&=&\dfrac{\lambda}{\pi}\int_{\R} e^{itx} \bigl(1+\beta\,\sgn(x)\bigr)\, \dfrac{1-\cos(hx)}{x^2} \dd x,\quad t\geqslant 0.
\end{eqnarray*}
Thus we may write
\begin{eqnarray*}
	\psi(t,h)=\int_{\R}e^{itx}\bigl(1-\cos(hx)\bigr)\tfrac{1+x^2}{x^2}\dd G(x),\quad t\geqslant 0,\,\,h>0,
\end{eqnarray*}
with 
\begin{eqnarray*}
	G(x):=\dfrac{\lambda}{\pi}\int_{-\infty}^{x} \dfrac{1+\beta\,\sgn(u)}{1+u^2}\,\dd u,\quad x\in\R.
\end{eqnarray*}
In the latter integral, the integrand function is non-negative  (recall that $\beta\in[-1,1]$) and it belongs to $L_1(\R)$. So $G\in\BV$ and it is non-decreasing on $\R$. According to Theorem \ref{th_crit_psi} with Corollary \ref{co_th_crit_psi} and Theorem~\ref{th_InfDiv}, the function $f$ defined by formula \eqref{def_f_Cauchy} is the characteristic for infinitely divisible law with the spectral components $\gamma$ (it equals $\Arg f(1)$) and  $G$.\quad $\Box$\\

\textbf{Example 4.} Here we illustrate application of Theorem \ref{th_suff_phi} with approximation. We consider the characteristic function of the Linnik distribution with $\alpha=1$ (see \cite{Linnik}):
\begin{eqnarray*}
	f(t)=\dfrac{1}{1+|t|},\quad t\in\R.
\end{eqnarray*}
It can be easily checked by P\'olya's condition that $f$ is indeed a characteristic function of absolutely continuous probability law. The density function, say $p$, admits the following representation (see \cite{Lukacs} p. 83):
\begin{eqnarray*}
	p(x)=\dfrac{1}{\pi}\lim_{R\to\infty}\int_0^R \cos(tx) f(t)\dd t =\dfrac{1}{\pi}\lim_{R\to\infty}\int_0^R  \dfrac{\cos(tx)}{1+t}\dd t,\quad x\in\R.
\end{eqnarray*}

We fix $h>0$ and consider the function $\varphi_+$ (see formula \eqref{def_phi}):
\begin{eqnarray*}
\varphi_+(t,h)=\dfrac{f(t)^2}{f(t-h)f(t+h)}-1&=& \dfrac{(1+|t-h|)(1+|t+h|)}{(1+|t|)^2}-1\\&=&\dfrac{|t-h|+|t+h|-2|t|+|t-h|\cdot|t+h|-|t|^2}{(1+|t|)^2},\quad t\in\R.
\end{eqnarray*}
Observe that $\varphi_+(t,h)=\tfrac{2h-2t+h^2-2t^2}{(1+t)^2}$ for any $t\in[0,h)$, and $\varphi_+(t,h)=-\tfrac{h^2}{(1+t)^2}$ for any $t\geqslant h$. So we rewrite $\varphi_+$ in the following form
\begin{eqnarray}\label{def_hatphi}
	\varphi_+(t,h)=2\hat\varphi_+(t,h)-h^2f(t)^2,\quad t\geqslant 0, 
\end{eqnarray}
where
\begin{eqnarray*}
	\hat\varphi_+(t,h)=
	\begin{cases}
		\dfrac{h-t+h^2-t^2}{(1+t)^2},& t<h,\\
		0,& t\geqslant h.
	\end{cases}
\end{eqnarray*}
According to Theorem \ref{th_suff_phi}, we have to propose a good approximation of the function $\varphi_+(\,\cdot\,,h)$ by some function $I(\,\cdot\,,h)$ that admits the suitable integral representation \eqref{def_I}. So we approximate  $\hat\varphi_+(\,\cdot\,,h)$ by the function $\Delta_h$ that has representation of such type according to formula \eqref{eq_Deltah}, and we introduce the corresponding approximation error $\rho_1(t,h):= \hat{\varphi}_+(t,h)-\Delta_h(t)$, $t\geqslant 0$. Next, recall that $t\mapsto f(t)^2$ is the characteristic function for $p^{*2}$. So the component $h^2 f(t)^2$ in \eqref{def_hatphi}, of course, admits the integral representation
\begin{eqnarray*}
	h^2 f(t)^2=h^2\int_{\R} e^{itx} p^{*2}(x) \dd x= 2\int_{\R} e^{itx} \dfrac{h^2}{2}\, p^{*2}(x) \dd x,\quad t\in\R,
\end{eqnarray*}
but it isn't of the form \eqref{def_I}. Therefore we have to replace $\tfrac{h^2}{2}$ with $\tfrac{1-\cos(hx)}{x^2}$ in the integral, i.e. we approximate $h^2 f(t)^2$ by
\begin{eqnarray*}
	2\int_{\R} e^{itx}\,\tfrac{1-\cos(hx)}{x^2}\,p^{*2}(x)\dd x =2\pi\int_{\R} e^{itx}D_h(x)\,p^{*2}(x)\dd x,\quad t\in\R,
\end{eqnarray*}
where $D_h$ is defined by formula \eqref{def_Dh}. Observe that $\Delta_h$, $f^2$ and $\Delta_h*f^2$ are from $L_1(\R)$ and we may write
\begin{eqnarray*}
	D_h(x)\,p^{*2}(x)=\dfrac{1}{2\pi}\int_{\R} e^{-itx} \Delta_h(t)\dd t \cdot  \dfrac{1}{2\pi}\int_{\R} e^{-itx} f(t)^2\dd t=\dfrac{1}{(2\pi)^2} \int_{\R} e^{-itx} (\Delta_h*f^2)(t)\dd t,\quad x\in\R.
\end{eqnarray*}
Here the function $D_h\cdot p^{*2}$ is from $L_1(\R)$ too and, inverting the Fourier transform, we get
\begin{eqnarray*}
	2\pi\int_{\R} e^{itx} D_h(x)\,p^{*2}(x)\dd x=(\Delta_h* f^2)(t),\quad t\in\R.
\end{eqnarray*}
Thus we may write the approximation error in the form: $\rho_2(t,h):=h^2 f(t)^2- (\Delta_h*f^2)(t)$, $t\geqslant 0$.

We now represent $\varphi_+$ as in \eqref{th_suff_phi_bkh}, i.e.
\begin{eqnarray*}
	\varphi_+(t,h)=2\hat\varphi_+(t,h)-h^2f(t)^2=2I(t,h)+\rho(t,h),\quad t\geqslant 0, 
\end{eqnarray*} 
where $I(t,h):= \Delta_h(t)-\tfrac{1}{2}\, (\Delta_h*f^2)(t)$, and
\begin{eqnarray*}
	\rho(t,h):=  2\hat\varphi_+(t,h)-2\Delta_h(t)-h^2 f(t)^2+(\Delta_h*f^2)(t)=2\rho_1(t,h)-\rho_2(t,h).
\end{eqnarray*}

Let us consider the function $I$ using the definition of $D_h$:
\begin{eqnarray*}
	I(t,h)&=&\int_{\R} e^{itx} D_h(x) \dd x -\pi\int_{\R} e^{itx}D_h(x)\,p^{*2}(x)\dd x\\
	&=&\int_{\R} e^{itx} D_h(x)(1-\pi p^{*2}(x)) \dd x
	=\int_{\R} e^{itx}\dfrac{1-\cos(hx)}{ x^2}\,\biggl(\dfrac{1}{\pi}-p^{*2}(x)\biggr) \dd x,\quad t\geqslant 0.
\end{eqnarray*}
Thus  we come to the needed representation for any $t\geqslant 0$ and $h>0$:
\begin{eqnarray*}
	I(t,h)=\int_{\R} e^{itx} (1-\cos(hx))\,\tfrac{1+x^2}{x^2}\,\dd G(x),
\end{eqnarray*}
with
\begin{eqnarray*}
	G(x):=\int_{-\infty}^{x}\biggl(\dfrac{1}{\pi}-p^{*2}(u)\biggr)\, \dfrac{\dd u}{(1+u^2)} ,\quad x\in\R,
\end{eqnarray*}
where the integrand function belongs to $L_1(\R)$ and it is always  non-negative. The latter fact is true, because $p^{*2}(x)=\tfrac{1}{2\pi}\int_{\R} e^{-itx} f(t)^2 \dd t$, $x\in\R$, and, consequently,
\begin{eqnarray*}
	p^{*2}(x)\leqslant\dfrac{1}{2\pi}\int_{\R} f(t)^2 \dd t=\dfrac{1}{2\pi}\int_{\R} \dfrac{\dd t}{(1+|t|)^2} = \dfrac{1}{\pi}\int_{0}^\infty  \dfrac{\dd t}{(1+t)^2} =\dfrac{1}{\pi},\quad x\in\R.
\end{eqnarray*}
Thus $G\in\BV$ and it is non-decreasing on $\R$.

We will show that the function $\rho$ satisfies condition \eqref{th_suff_phi_b} of Theorem \ref{th_suff_phi}. We fix any sequence $(h_l)_{l\in\N}$ as in the statement of Theorem \ref{th_suff_phi} and we first observe that
\begin{eqnarray*}
	|\rho(kh_l,h_l)|= |2\rho_1(kh_l,h_l)-\rho_2(kh_l,h_l)|\leqslant 2|\rho_1(kh_l,h_l)|+|\rho_2(kh_l,h_l)|,\quad k\in\N_0.
\end{eqnarray*}
Here and below $l\in\N$ is choosen arbitrarily. So $\rho_1(kh_l,h_l)=\hat{\varphi}_+(kh_l,h_l)-\Delta_{h_l}(kh_l)$ equals $(h_l^2+h_l)-h_l=h^2_l$ for $k=0$ and it equals $0$ else. For the function $\rho_2$  we preliminarily observe that for any $t\geqslant0$ and $h>0$
\begin{eqnarray*}
	\rho_2(t,h)=\int_{-h}^{h} \Delta_{h}(s)\dd s\cdot f(t)^2-\int_{-h}^{h}f^2(t-s) \Delta_{h}(s)\dd s=\int_{-h}^{h}\bigl(f^2(t)-f^2(t-s)\bigr) \Delta_{h}(s)\dd s.
\end{eqnarray*}
Here
\begin{eqnarray*}
	f^2(t)-f^2(t-s)&=& \dfrac{1}{(1+|t|)^2}- \dfrac{1}{(1+|t-s|)^2}\\
	&=&\dfrac{(1+|t-s|)^2-(1+|t|)^2}{(1+|t|)^2(1+|t-s|)^2}=\dfrac{(|t-s|-|t|)(2+|t-s|+|t|)}{(1+|t|)^2(1+|t-s|)^2},
\end{eqnarray*}
and we have the estimate
\begin{eqnarray*}
	|f^2(t)-f^2(t-s)|\leqslant \dfrac{|s|\cdot 2(1+|t|+|t-s|)}{(1+|t|)^2(1+|t-s|)^2}\leqslant \dfrac{2|s|}{(1+|t|)(1+|t-s|)}\leqslant 2|s|.
\end{eqnarray*}
Then for any $t\geqslant 0$ and $h>0$ we get
\begin{eqnarray*}
	|\rho_2(t,h)|\leqslant \int_{-h}^{h}|f^2(t)-f^2(t-s)| \Delta_{h}(s)\dd s\leqslant  \int_{-h}^{h}2|s| \Delta_{h}(s)\dd s\leqslant 2h  \int_{-h}^{h} \Delta_{h}(s)\dd s=2h^3.
\end{eqnarray*}
In particular, $|\rho_2(kh_l,h_l)|\leqslant 2h_l^3$. We now turn to checking \eqref{th_suff_phi_b} for any $t>0$ using the remarks about $\rho_1$ and $\rho_2$:
\begin{eqnarray*}
	\sum_{k=0}^{n_{t,l}-1}\bigl(n_{t,l}-k\bigr)|\rho(kh_l,h_l)|&\leqslant& 2\sum_{k=0}^{n_{t,l}-1}\bigl(n_{t,l}-k\bigr) |\rho_1(kh_l,h_l)|+\sum_{k=0}^{n_{t,l}-1}\bigl(n_{t,l}-k\bigr)|\rho_2(kh_l,h_l)|\\
	&\leqslant& 2 n_{t,l} h_l^2+\sum_{k=0}^{n_{t,l}-1}\bigl(n_{t,l}-k\bigr)| 2h_l^3\leqslant  2 n_{t,l} h_l^2+2n_{t,l}^2 h_l^3.
\end{eqnarray*}
Since $n_{t,l}= \lfloor t/ h_l\rfloor\leqslant t/ h_l$, we have
\begin{eqnarray*}
	\sum_{k=0}^{n_{t,l}-1}\bigl(n_{t,l}-k\bigr)|\rho(kh_l,h_l)|\leqslant  2 t h_l+2t^2h_l,
\end{eqnarray*}
where the sum in the right-hand side tends to zero as $l\to\infty$ due to $h_l\to 0$ by the assumption. Thus \eqref{th_suff_phi_b} holds.

We still have to check condition \eqref{th_suff_phi_phi2} for the function $\varphi_+$. According to the comment before formula \eqref{def_hatphi}, we have $\varphi_+(0,h_l)=2h_l+h_l^2$, and $\varphi_+(kh_l,h_l)=-h_l^2/(1+kh_l)^2$ for any $k\in\N$, i.e.  $|\varphi_+(kh_l,h_l)|\leqslant h_l^2$, $k\in\N$. Then for any $t>0$ we have
\begin{eqnarray*}
	\sum_{k=0}^{n_{t,l}-1}\bigl(n_{t,l}-k\bigr)\bigl|\varphi_{+}(kh_l,h_l)\bigr|^2 \leqslant n_{t,l} \bigl(2h_l+h_l^2\bigr)^2+\sum_{k=1}^{n_{t,l}-1} \bigl(n_{t,l}-k\bigr)h_l^4\leqslant n_{t,l} h_l^2(2+h_l)^2+n_{t,l}^2 h_l^4.
\end{eqnarray*}
Using inequality $n_{t,l}\leqslant t/ h_l$, we obtain the estimate
\begin{eqnarray*}
	\sum_{k=0}^{n_{t,l}-1}\bigl(n_{t,l}-k\bigr)\bigl|\varphi_{+}(kh_l,h_l)\bigr|^2 \leqslant  t h_l(2+h_l)^2+t^2 h_l^2,
\end{eqnarray*}
where the right-hand side obviously tends to zero as $l\to\infty$ due to $h_l\to 0$. This implies \eqref{th_suff_phi_phi2} for $\varphi_+$.

Thus we have checked all conditions of Theorem \ref{th_suff_phi}, and we conclude that the distribution function of the Linnik law with characteristic function $f$ is rational-infinitely divisible, i.e. $F\in\Q$. Moreover, since $G$ is non-decreasing (as we have already showed above), $F$ is purely infinitely divisible, i.e. $F\in\ID$, with the spectral pair $\gamma=\Arg f(1)=0$ and $G$.

It should be noted that the infinite divisibility of the Linnik distributions (with characteristic functions $f(t)=1/(1+|t|^\alpha)$, $t\in\R$, with $\alpha\in(0,2)$) is a known fact. It was showed in the papers \cite{Devro} and \cite{KotzOstr}, where the authors used  other non-trivial approaches.
\quad $\Box$\\

\section{Proofs}

\textbf{Proof of Theorem \ref{th_crit_psi}.} We first show that for any $n\in\N_0$ and $h\in\R$ the quantity $\Ln f(nh)-in\Arg f(h)$ is expressed through the values $\Delta^2_h\Ln f(jh)$, $j\in\N_0$. We fix such $n$ and $h$. Observe that for any  $k\in\N$ 
\begin{eqnarray*}
	\sum_{j=1}^{k-1} \Delta^2_h\Ln f(jh)&=&\sum_{j=1}^{k-1} \bigl(\Ln f((j-1)h)+\Ln f((j+1)h)-2\Ln f(jh)\bigr)\\
	&=&\sum_{j=1}^{k-1} \bigl(\Ln f((j-1)h)-\Ln f(jh)\bigr)+\sum_{j=1}^{k-1} \bigl(\Ln f((j+1)h)-\Ln f(jh)\bigr)\\
	&=&\Ln f(0)-\Ln f((k-1)h)+\Ln f(kh)-\Ln f(h)\\
	&=&\Ln f(kh)-\Ln f((k-1)h)-\Ln f(h),
\end{eqnarray*}
and
\begin{eqnarray*}
	\sum_{j=0}^{k-1} \Delta^2_h\Ln f(jh)&=&  \sum_{j=1}^{k-1} \Delta^2_h\Ln f(jh)+  \Delta^2_h\Ln f(0)\\
	&=&\bigl(\Ln f(kh)-\Ln f((k-1)h)-\Ln f(h)\bigr)+ \Ln f(-h)+\Ln f(h)-2\Ln f(0)
	\\
	&=&\Ln f(kh)-\Ln f((k-1)h)+\Ln f(-h).
\end{eqnarray*}
Hence for any $k\in\N$
\begin{eqnarray*}
\sum_{j=1}^{k-1} \Delta^2_h\Ln f(jh)+\tfrac{1}{2}\Delta^2_h\Ln f(0)
	&=&\Ln f(kh)-\Ln f((k-1)h)-\tfrac{1}{2}\bigl(\Ln f(h)-\Ln f(-h)\bigr).
\end{eqnarray*}
Let us consider
\begin{eqnarray*}
\sum_{k=1}^{n} \biggl(\sum_{j=1}^{k-1} \Delta^2_h\Ln f(jh)+\tfrac{1}{2}\Delta^2_h\Ln f(0)\biggr)
	&=& \sum_{k=1}^{n} \bigl(\Ln f(kh)-\Ln f((k-1)h)\bigr)-\sum_{k=1}^{n}\tfrac{1}{2}\bigl(\Ln f(h)-\Ln f(-h)\bigr)\\
&=&\Ln f(nh)-\Ln f(0)-n\cdot\tfrac{1}{2}\bigl(\Ln f(h)-\Ln f(-h)\bigr)\\
&=&\Ln f(nh)-n\cdot\tfrac{1}{2}\bigl(\Ln f(h)-\Ln f(-h)\bigr).
\end{eqnarray*}
Since $\ln|f(\cdot)|$ and $\Arg f(\cdot)$ are respectively even and odd functions, we have 
\begin{eqnarray*}
	\Ln f(h)-\Ln f(-h)&=&\ln |f(h)|+i \Arg f(h) - \ln |f(-h)|-i\Arg f(-h)\\
	&=&\ln |f(h)|+i \Arg f(h) - \ln |f(h)|+i\Arg f(h)\\
	&=&2i \Arg f(h).
\end{eqnarray*}
Thus we obtain
\begin{eqnarray*}
	\Ln f(nh)-in\Arg f(h)&=&\sum_{k=1}^{n} \Bigl(\sum_{j=1}^{k-1} \Delta^2_h\Ln f(jh)+\tfrac{1}{2}\Delta^2_h\Ln f(0)\Bigr),\quad n\in\N_0, \,\, h\in\R,
\end{eqnarray*}
as required. By the way, we observe that for any $n\in\N_0$, $h\in\R$, and $m\in\N$:
\begin{eqnarray*}
	\Ln f(nh)-i\tfrac{n}{m}\Arg f(mh)&=&\Ln f(nh)-i\tfrac{n}{m}\Imagpart\{\Ln f(mh)\}\\
	&=&\Ln f(nh)-in \Arg f(h)-i\tfrac{n}{m}\Imagpart\bigl\{\Ln f(mh)-im \Arg f(h)\bigr\}\\
	&=& \sum_{k=1}^{n} \Bigl(\sum_{j=1}^{k-1} \Delta^2_h\Ln f(jh)+\tfrac{1}{2}\Delta^2_h\Ln f(0)\Bigr)\\
	&&{}-i\tfrac{n}{m}\Imagpart\biggl\{ \sum_{k=1}^{m} \Bigl(\sum_{j=1}^{k-1} \Delta^2_h\Ln f(jh)+\tfrac{1}{2}\Delta^2_h\Ln f(0)\Bigr)\biggr\}.
\end{eqnarray*}
Thus $\Ln f(nh)-i\tfrac{n}{m}\Arg f(mh)$ is also expressed through the values $\Delta^2_h\Ln f(jh)$, $j\in\N_0$.

We rewrite $\Ln f(nh)-i\tfrac{n}{m}\Arg f(mh)$ in terms of $\psi(jh,h)$ according to \eqref{def_Delta2Ln}:
\begin{eqnarray*}
	\Ln f(nh)-i\tfrac{n}{m}\Arg f(mh)&=&-2\sum_{k=1}^{n} \Bigl(\sum_{j=1}^{k-1} \psi(jh,h)+\tfrac{1}{2}\psi(0,h)\Bigr)\\
	&&{}+2i\tfrac{n}{m}\Imagpart \biggl\{ \sum_{k=1}^{m} \Bigl(\sum_{j=1}^{k-1} \psi(jh,h)+\tfrac{1}{2}\psi(0,h)\Bigr)\biggr\}.
\end{eqnarray*}
Let us turn to the functions $I$ and $\delta$ from \eqref{def_Idelta}. We have
\begin{eqnarray*}
	I(kh,h)= \int_{\R}e^{ikhx}\bigl(1-\cos(hx)\bigr)\tfrac{1+x^2}{x^2}\dd G(x)=\int_{\R} e^{ikhx}  \dd R_h(x)\quad\text{for all}\quad k\in\N_0,\,\, h\in\R,
\end{eqnarray*}
where
\begin{eqnarray}\label{def_Rh}
	R_h(x)=  \int_{(-\infty,x]}\bigl(1-\cos(hu)\bigr)\tfrac{1+u^2}{u^2}\dd G(u),\quad x\in\R.
\end{eqnarray}
So $\psi(kh,h)=I(kh,h)+\delta(kh,h)$, $k\in\N_0$, $h\in\R$. Let us represent
\begin{eqnarray}\label{repr_Lnfnmh}
	\Ln f(nh)-i\tfrac{n}{m}\Arg f(mh)&=&-2\sum_{k=1}^{n} \Bigl(\sum_{j=1}^{k-1} I(jh,h)+\tfrac{1}{2}I(0,h)\Bigr)\nonumber\\
	&&{}+2i\tfrac{n}{m}\Imagpart\biggl\{ \sum_{k=1}^{m} \Bigl(\sum_{j=1}^{k-1} I(jh,h)+\tfrac{1}{2}I(0,h)\Bigr)\biggr\}-2\Delta(n,m,h)
\end{eqnarray}
for any $n\in\N_0$, $m\in\N$, $h\in\R$, where $\Delta(n,m,h)$ admits the representation:
\begin{eqnarray*}
	\Delta(n,m,h)=\sum_{k=1}^{n} \Bigl(\sum_{j=1}^{k-1} \delta(jh,h)+\tfrac{1}{2}\delta(0,h)\Bigr)-i\tfrac{n}{m}\Imagpart\biggl\{ \sum_{k=1}^{m} \Bigl(\sum_{j=1}^{k-1}\delta(jh,h)+\tfrac{1}{2}\delta(0,h)\Bigr)\biggr\}.
\end{eqnarray*}
For this we have the estimates:
\begin{eqnarray}\label{ineq_Deltanmh}
	|\Delta(n,m,h)|&\leqslant& \sum_{k=1}^{n} \Bigl(\sum_{j=1}^{k-1} |\delta(jh,h)|+\tfrac{1}{2}|\delta(0,h)|\Bigr)+\tfrac{n}{m} \sum_{k=1}^{m} \Bigl(\sum_{j=1}^{k-1} |\delta(jh,h)|+\tfrac{1}{2}|\delta(0,h)|\Bigr)\nonumber\\
	&\leqslant& \sum_{k=1}^{n} \sum_{j=0}^{k-1} |\delta(jh,h)|+\tfrac{n}{m} \sum_{k=1}^{m} \sum_{j=0}^{k-1} |\delta(jh,h)|\nonumber\\
	&=& \sum_{k=0}^{n-1} (n-k) |\delta(kh,h)|+\tfrac{n}{m} \sum_{k=0}^{m-1} (m-k) |\delta(kh,h)|.
\end{eqnarray}

We now fix any $n\in\N_0$, $m\in\N$, $h>0$, and we consider the sums
\begin{eqnarray*}
	\sum_{k=1}^{n} \Bigl(\sum_{j=1}^{k-1} I(jh,h)+\tfrac{1}{2}I(0,h)\Bigr)=\sum_{k=1}^{n} \Bigl(\sum_{j=0}^{k-1} I(jh,h)-\tfrac{1}{2}I(0,h)\Bigr)=\sum_{k=1}^{n} \sum_{j=0}^{k-1} I(jh,h)-\tfrac{n}{2}I(0,h).
\end{eqnarray*}
For any $k\in\N$
\begin{eqnarray*}
	\sum_{j=0}^{k-1} I(jh,h)=\sum_{j=0}^{k-1} \int_{\R} e^{ijh x}  \dd R_h(x)
	= \int_{\R} \sum_{j=0}^{k-1}e^{ijhx}  \dd R_h(x)=\int_{\R} \dfrac{e^{ikhx}-1}{e^{ihx}-1} \dd R_h(x).
\end{eqnarray*}
Next, we write
\begin{eqnarray*}
	\sum_{k=1}^{n} \sum_{j=0}^{k-1} I(jh,h)=\sum_{k=1}^{n} \int_{\R} \dfrac{e^{ikhx}-1}{e^{ihx}-1} \dd R_h(x)= \int_{\R} \sum_{k=1}^{n}\dfrac{e^{ikhx}-1}{e^{ihx}-1} \dd R_h(x)=\int_{\R} \dfrac{\sum_{k=1}^{n}e^{ikhx}-n}{e^{ihx}-1} \dd R_h(x).
\end{eqnarray*}
Observe that
\begin{eqnarray*}
	\sum_{k=1}^{n}e^{ikhx}=\dfrac{e^{i(n+1)hx}-1}{e^{ihx}-1} -1.
\end{eqnarray*}
Hence
\begin{eqnarray*}
	\dfrac{\sum_{k=1}^{n}e^{ikhx}-n}{e^{ihx}-1} = \dfrac{e^{i(n+1)hx}-1}{(e^{ihx}-1)^2}- \dfrac{n+1}{e^{ihx}-1}=\dfrac{e^{i(n+1)hx}-1-(n+1)(e^{ihx}-1)}{(e^{ihx}-1)^2}.
\end{eqnarray*}
Thus
\begin{eqnarray*}
	\sum_{k=1}^{n} \sum_{j=0}^{k-1} I(jh,h)= \int_{\R}  \biggl( \dfrac{e^{i(n+1)hx}-1-(n+1)(e^{ihx}-1)}{(e^{ihx}-1)^2}\biggr)  \dd R_h(x).
\end{eqnarray*}
Since
\begin{eqnarray*}
	I(0,h)=\int_{\R}\dd R_h(x),
\end{eqnarray*}
we have
\begin{eqnarray*}
	\sum_{k=1}^{n} \sum_{j=0}^{k-1} I(jh,h)-\tfrac{n}{2}I(0,h)
	=\int_{\R}  \biggl( \dfrac{e^{i(n+1)hx}-1-(n+1)(e^{ihx}-1)}{(e^{ihx}-1)^2} -\dfrac{n}{2}\biggr)  \dd R_h(x).
\end{eqnarray*}
Let us simplify the expression in the integral. We first observe that
\begin{eqnarray*}
	|e^{ihx}-1|^2=(1-\cos (hx))^2+(\sin (hx))^2=(1-\cos (hx))^2+1-(\cos (hx))^2=2(1-\cos (hx)).
\end{eqnarray*}
Next, we write
\begin{eqnarray*}
	|e^{ihx}-1|^2=(e^{ih}-1)(e^{-ih}-1)= e^{-ihx}(e^{ihx}-1)(1-e^{ihx})=-e^{-ihx}(e^{ihx}-1)^2.
\end{eqnarray*}
Therefore
\begin{eqnarray*}
	(e^{ihx}-1)^2= -e^{ihx}\cdot 2(1-\cos (hx)).
\end{eqnarray*}
So we write
\begin{eqnarray*}
	\dfrac{e^{i(n+1)hx}-1-(n+1)(e^{ihx}-1)}{(e^{ihx}-1)^2}
	&=&-\dfrac{\bigl(e^{i(n+1)hx}-1-(n+1)(e^{ihx}-1)\bigr)e^{-ihx}}{ 2(1-\cos (hx))}\\
	&=&-\dfrac{e^{inhx}-e^{-ihx}-(n+1)(1-e^{-ihx})}{ 2(1-\cos (hx))}\\
	&=&-\dfrac{e^{inhx}-1-n(1-e^{-ihx})}{ 2(1-\cos (hx))}.
\end{eqnarray*}
Here $1-e^{-ihx}=1-\cos(hx)+i\sin(hx)$. Then
\begin{eqnarray*}
	-\dfrac{e^{inhx}-1-n(1-e^{-ihx})}{ 2(1-\cos (hx))}= -\dfrac{e^{inhx}-1-in\sin(hx)}{ 2(1-\cos (hx))}+\dfrac{n}{2},
\end{eqnarray*}
and
\begin{eqnarray*}
	\dfrac{e^{i(n+1)hx}-1-(n+1)(e^{ihx}-1)}{(e^{ihx}-1)^2}-\dfrac{n}{2}
	=-\dfrac{e^{inhx}-1-in\sin(hx)}{ 2(1-\cos (hx))}.
\end{eqnarray*}
Thus
\begin{eqnarray*}
	\sum_{k=1}^{n} \sum_{j=0}^{k-1} I(jh,h)-\tfrac{n}{2}I(0,h)
	=\int_{\R}  \biggl( -\dfrac{e^{inhx}-1-in\sin(hx)}{ 2(1-\cos (hx))}\biggr)  \dd R_h(x).
\end{eqnarray*}
Returning  to \eqref{repr_Lnfnmh}, we analogously write
\begin{eqnarray*}
	\sum_{k=1}^{m} \sum_{j=0}^{k-1} I(jh,h)-\tfrac{m}{2}I(0,h)
	=\int_{\R}  \biggl( -\dfrac{e^{imhx}-1-im\sin(hx)}{ 2(1-\cos (hx))}\biggr)  \dd R_h(x),
\end{eqnarray*}
and also
\begin{eqnarray*}
	\Imagpart\biggl\{\sum_{k=1}^{m} \sum_{j=0}^{k-1} I(jh,h)-\tfrac{m}{2}I(0,h)\biggr\}
	=\int_{\R}  \biggl(-\dfrac{\sin(mhx)-m\sin(hx)}{2(1-\cos (hx))}\biggr)  \dd R_h(x).
\end{eqnarray*}
Hence, due to \eqref{repr_Lnfnmh}, we have
\begin{eqnarray*}
	\Ln f(nh)-i\tfrac{n}{m}\Arg f(mh)&=&\int_{\R}   \dfrac{e^{inhx}-1-in\sin(hx)}{ 1-\cos (hx)} \dd R_h(x)\\
	&&{}\quad-i\tfrac{n}{m}\cdot \int_{\R}   \dfrac{\sin(mhx)-m\sin(hx)}{ 1-\cos (hx)} \dd R_h(x)-2\Delta(n,m,h)\\
	&=& \int_{\R}  \dfrac{e^{inhx}-1-i\tfrac{n}{m}\sin(mhx)}{ 1-\cos (hx)} \dd R_h(x)-2\Delta(n,m,h).
\end{eqnarray*}
According to \eqref{def_Rh}, we obtain the representation
\begin{eqnarray}\label{repr_Ln_int_Delta}
	\Ln f(nh)-i\tfrac{n}{m}\Arg f(mh)= \int_{\R}  \bigl( e^{inhx}-1-i\tfrac{n}{m}\sin(mhx)\bigr)\tfrac{1+x^2}{x^2}\dd G(x)-2\Delta(n,m,h)
\end{eqnarray}
for any $n\in\N_0$, $m\in\N$, and $h>0$.

Let  us  choose the sequence $(h_l)_{l\in\N}$ and integers $n_{t,l}$ as in the statement of the theorem, i.e. $h_l\to 0$ as $l\to\infty$, $n_{t,l}\defeq \lfloor t/ h_l\rfloor\in \N_0$, $t>0$, $l\in\N$, and also \eqref{th_crit_psi_cond_a} holds. Without loss of generality we can assume that $h_l\leqslant \tfrac{1}{2}$ for all $l\in\N$. Let us define $m_l\defeq n_{1,l}=\lfloor 1/ h_l\rfloor\in\N$, $l\in\N$.    We fix $t>0$. Observe that  $0\leqslant t-n_{t,l} h_l<h_l$ and $0\leqslant 1-m_{l} h_l<h_l$ for any $l\in\N$. Since $n_{t,l} h_l\to t$ and $m_lh_l\to 1$ as $l\to\infty$, we have
\begin{eqnarray*}
	\Ln f(t)-it\Arg f(1)&=& \lim_{l\to\infty}\biggl( \Ln f(n_{t,l} h_l)-in_{t,l} h_l\Arg f(1)\biggr)\\
	&=& \lim_{l\to\infty}\biggl( \Ln f(n_{t,l} h_l)-in_{t,l} h_l\,\dfrac{\Arg f(m_l h_l)}{m_l h_l}\biggr)\\
	&=& \lim_{l\to\infty}\Bigl( \Ln f(n_{t,l}h_l)-i \tfrac{n_{t,l}}{m_l}\,\Arg f(m_lh_l)\Bigr).
\end{eqnarray*}
Using \eqref{repr_Ln_int_Delta}, we write	
\begin{eqnarray*}
	\Ln f(t)-it\Arg f(1)	
	=\lim_{l\to\infty} \Biggl(\,\int_{\R}\bigl( e^{in_{t,l} h_l x}-1-\tfrac{in_l}{m_l}\sin(m_l h_l x)\bigr)\tfrac{1+x^2}{x^2}\dd G(x)-2\Delta(n_{t,l},m_l,h_l)\Biggr).
\end{eqnarray*}

We first consider the sequence $\Delta(n_{t,l},m_l,h_l)$, $l\in\N$. Applying inequality \eqref{ineq_Deltanmh}, we get
\begin{eqnarray*}
	|\Delta(n_{t,l},m_l,h_l)|\leqslant\sum_{k=0}^{n_{t,l}-1} (n_{t,l}-k) |\delta(kh_l,h_l)|+\tfrac{n_{t,l}}{m_l} \sum_{k=0}^{m_l-1} (m_l-k) |\delta(kh_l,h_l)|.
\end{eqnarray*}
Due to \eqref{th_crit_psi_cond_a},  the first sum goes to 0 as $l\to\infty$. The second sum goes to 0  as $l\to\infty$ too, because  $m_l= n_{1,l}$ and we again apply \eqref{th_crit_psi_cond_a} with $t=1$. Also $\tfrac{n_{t,l}}{m_l} =\tfrac{n_{t,l}h_l}{m_l h_l}\to t$ as $l\to\infty$. Thus we showed that $\Delta(n_{t,l},m_l,h_l)\to 0$ as $l\to\infty$.

We now have
\begin{eqnarray}\label{conc_Lnlim}
	\Ln f(t)-it\gamma
	= \lim_{l\to\infty} \int_{\R}\bigl( e^{in_{t,l} h_l x}-1-\tfrac{in_l}{m_l}\sin(m_l h_l x)\bigr)\tfrac{1+x^2}{x^2}\dd G(x),
\end{eqnarray}
where we set $\gamma:=\Arg f(1)$.
It is easily seen that
\begin{eqnarray}\label{conc_integrand}
\lim_{l\to\infty}\Bigl[\bigl( e^{in_{t,l} h_l x}-1-\tfrac{in_l}{m_l}\sin(m_l h_l x)\bigr)\tfrac{1+x^2}{x^2}\Bigr]=\bigl( e^{i t x}-1-it\sin(x)\bigr)\tfrac{1+x^2}{x^2}\quad\text{for any}\quad x\in\R.
\end{eqnarray}

We will dominate the integrand in \eqref{conc_Lnlim}. Let $r$ be a positive real number. For any $x\in\R$ satisfying  $|x|> r$ we have the estimate
\begin{eqnarray*}
	\bigl| e^{in_{t,l} h_l x}-1-\tfrac{in_l}{m_l}\sin(m_l h_l x)\bigr|\tfrac{1+x^2}{x^2}\leqslant  \Bigl( \bigl|e^{in_{t,l} h_l x}\bigr|+1+\bigl|\tfrac{in_l}{m_l}\sin(m_l h_l x)\bigr|\Bigr)\tfrac{1+x^2}{x^2}\leqslant  \bigl(2+\tfrac{n_l}{m_l}\bigr)\tfrac{1+r^2}{r^2}.
\end{eqnarray*}
Since $n_lh_l\leqslant t$ and $m_l h_l> 1-h_l\geqslant \tfrac{1}{2}$, we have $\tfrac{n_l}{m_l}= \tfrac{n_lh_l}{m_lh_l}\leqslant 2t$. Thus 
\begin{eqnarray}\label{ineq_estim_great_r}
	\bigl| e^{in_{t,l} h_l x}-1-\tfrac{in_l}{m_l}\sin(m_l h_l x)\bigr|\tfrac{1+x^2}{x^2}\leqslant    2(1+t)\tfrac{1+r^2}{r^2},\quad |x|>r,\quad l\in\N.
\end{eqnarray}
For the case $|x|\leqslant r$ we will use the representation
\begin{eqnarray*}
	e^{in_{t,l} h_l x}-1-\tfrac{in_l}{m_l}\sin(m_l h_l x)= \bigl(e^{in_{t,l} h_l x}-1-in_{t,l} h_l x\bigr) -\tfrac{in_l}{m_l}\bigl(\sin(m_l h_l x)- m_l h_l x\bigr).
\end{eqnarray*}
Then we have the estimate
\begin{eqnarray*}
	\bigl|e^{in_{t,l} h_l x}-1-\tfrac{in_l}{m_l}\sin(m_l h_l x)\bigr|\leqslant \bigl|e^{in_{t,l} h_l x}-1-in_{t,l} h_l x\bigr| +\tfrac{n_l}{m_l}\bigl|\sin(m_l h_l x)- m_l h_l x\bigr|.
\end{eqnarray*}
Applying the well-known inequalities $|e^{iy}-1-iy|\leqslant \tfrac{y^2}{2!}$ and $|\sin y - y|\leqslant \tfrac{|y|^3}{3!}$, $y\in\R$,  we get
\begin{eqnarray*}
	\bigl|e^{in_{t,l} h_l x}-1-in_{t,l} h_l x\bigr|\leqslant \dfrac{(n_{t,l} h_l x)^2}{2!},\qquad \bigl|\sin(m_l h_l x)- m_l h_l x\bigr|\leqslant \dfrac{|m_l h_l x|^3}{3!}.
\end{eqnarray*}
Therefore
\begin{eqnarray*}
	\bigl| e^{in_{t,l} h_l x}-1-\tfrac{in_l}{m_l}\sin(m_l h_l x)\bigr|\tfrac{1+x^2}{x^2}\leqslant  \Bigl( \tfrac{(n_{t,l} h_l x)^2}{2!}+ \tfrac{|m_l h_l x|^3}{3!}\Bigr)\tfrac{1+x^2}{x^2}\leqslant  \Bigl( \tfrac{(n_{t,l} h_l)^2}{2!}+ \tfrac{(m_l h_l )^3 |x|}{3!}\Bigr)(1+x^2).
\end{eqnarray*}
Since $|x|\leqslant r$, $n_{t,l} h_l\leqslant t$,  and $m_l h_l\leqslant 1$, we obtain
\begin{eqnarray}\label{ineq_estim_less_r}
	\bigl| e^{in_{t,l} h_l x}-1-\tfrac{in_l}{m_l}\sin(m_l h_l x)\bigr|\tfrac{1+x^2}{x^2}\leqslant  \Bigl( \tfrac{t^2}{2!}+ \tfrac{r}{3!}\Bigr)(1+r^2),\quad |x|\leqslant r, \quad l\in\N.
\end{eqnarray}
Thus \eqref{ineq_estim_great_r} and \eqref{ineq_estim_less_r} together yield the needed domination of the intagrand in \eqref{conc_Lnlim}:
\begin{eqnarray*}
	\bigl| e^{in_{t,l} h_l x}-1-\tfrac{in_l}{m_l}\sin(m_l h_l x)\bigr|\tfrac{1+x^2}{x^2}\leqslant 	C_{t,r} ,\quad x\in\R, \quad l\in\N,
\end{eqnarray*}
with the following constant majorant (for fixed $t$ and $r$):
\begin{eqnarray*}
	C_{t,r}\defeq \max\Bigl\{ \tfrac{2(1+t)}{r^2} ,\tfrac{t^2}{2!}+ \tfrac{r}{3!}\Bigr\}\cdot(1+r^2).
\end{eqnarray*}

We now return to \eqref{conc_Lnlim}. Due  to the Jordan decomposition, we decompose $G(x)=G^+(x)-G^-(x)$, $x\in\R$,  with some non-decreasing functions $G^+$ and $G^-$ from $\BV$. So we have
\begin{eqnarray*}
	\Ln f(t)-it\gamma	
	&=& \lim_{l\to\infty}\Biggl( \int_{\R}\bigl( e^{in_{t,l} h_l x}-1-\tfrac{in_l}{m_l}\sin(m_l h_l x)\bigr)\tfrac{1+x^2}{x^2}\dd G^+(x)\\
	&&{}\qquad\qquad-\int_{\R}\bigl( e^{in_{t,l} h_l x}-1-\tfrac{in_l}{m_l}\sin(m_l h_l x)\bigr)\tfrac{1+x^2}{x^2}\dd G^-(x)\Biggr).
\end{eqnarray*}
Due to the convergence \eqref{conc_integrand} and the proved domination, by Lebesgue's dominated convergence theorem, we obtain
\begin{eqnarray*}
	\Ln f(t)-it\gamma	
	&=&\int_{\R}\bigl( e^{i t x}-1-it\sin(x)\bigr)\tfrac{1+x^2}{x^2}\dd G^+(x)
	-\int_{\R}\bigl( e^{i t x}-1-it\sin(x)\bigr)\tfrac{1+x^2}{x^2}\dd G^-(x)\\
	&=&\int_{\R}\bigl( e^{i t x}-1-it\sin(x)\bigr)\tfrac{1+x^2}{x^2}\dd G(x),
\end{eqnarray*}
where $t>0$ was fixed arbitrarily.

Let us show that the obtained formula holds for $t\leqslant 0$ too. Since $\Ln f(0)=0$, it is easily seen that the formula is true for $t=0$.  If $t<0$  then $-t>0$ and we can write
\begin{eqnarray*}
	\Ln f(-t)-i(-t)\gamma	= \int_{\R}\bigl( e^{-i t x}-1-i(-t)\sin(x)\bigr)\tfrac{1+x^2}{x^2}\dd G(x).	
\end{eqnarray*}
Hence, according to the properties of the components $\ln |f(\cdot)|$ and $\Arg f(\cdot)$ of the logarithm $\Ln f(\cdot)$, we get
\begin{eqnarray*}
	\Ln f(t)-it\gamma	=\overline{\Ln f(-t)-i(-t)\gamma	}=\int_{\R}\bigl( e^{i t x}-1-it\sin(x)\bigr)\tfrac{1+x^2}{x^2}\dd G(x).
\end{eqnarray*} 
Thus we obtain that
\begin{eqnarray*}
	\Ln f(t)=it\gamma+\int_{\R}\bigl( e^{i t x}-1-it\sin(x)\bigr)\tfrac{1+x^2}{x^2}\dd G(x),\quad t\in\R,
\end{eqnarray*}
i.e. \eqref{repr_f} holds as required.

If we additionally know that $f$ is the characteristic function of some distribution function $F$, then $F$ is quasi-infinitely divisible and, consequently, $F\in\Q$ (see Introduction).\quad $\Box$\\

\textbf{Proof of Theorem \ref{th_suff_phi}.} According to the definitions of $\psi$ and $\varphi_{\pm}$ we have
\begin{eqnarray*}
	\psi(t,h)=\pm\tfrac{1}{2}\,\Ln\bigl(1+\varphi_{\pm}(t,h)\bigr),\quad t\in\R,\,\, h\in\R.
\end{eqnarray*}

Next, we define a function $\delta$ by the formula \eqref{def_Idelta} with a given $I$, i.e. $\delta(t, h)=\psi(t,h)-I(t,h)$, $t\geqslant 0$, $h>0$. On account of \eqref{th_suff_phi_bkh}, for any $t\geqslant 0$ and $h>0$ we accordingly have
\begin{eqnarray*}
	\delta(t,h)= \psi(t,h) \mp  \tfrac{1}{2}\bigl(\varphi_{\pm}(t,h)-\rho(t,h)\bigr)
	= \pm\tfrac{1}{2}\,\Ln\bigl(1+\varphi_{\pm}(t,h)\bigr) \mp  \tfrac{1}{2}\,\varphi_{\pm}(t,h)\pm \tfrac{1}{2} \rho(t,h).
\end{eqnarray*}
Then
\begin{eqnarray*}
	|\delta(t,h)|\leqslant \tfrac{1}{2}\bigl|\Ln\bigl(1+\varphi_{\pm}(t,h)\bigr)-  \varphi_{\pm}(t,h)\bigr|+ \tfrac{1}{2}|\rho(t,h)|,\quad t\geqslant 0,\,\, h>0.
\end{eqnarray*}
The first modulus less than $| \varphi_{\pm}(t,h)|^2$, when $| \varphi_{\pm}(t,h)|\leqslant \tfrac{1}{2}$. Indeed, under this inequality, we have
\begin{eqnarray*}
	\Ln\bigl(1+\varphi_{\pm}(t,h)\bigr)=\ln\bigl(1+\varphi_{\pm}(t,h)\bigr)=\sum_{k=1}^{\infty} \tfrac{(-1)^{k-1}}{k}\, \varphi_{\pm}(t,h)^k,
\end{eqnarray*}
where the function  $\ln(\cdot)$   returns as usual  the principal value of the logarithm.
Therefore
\begin{eqnarray*}
	\bigr|\Ln\bigl(1+\varphi_{\pm}(t,h)\bigr)- \varphi_{\pm}(t,h)\bigr|\leqslant\dfrac{1}{2} \sum_{k=2}^{\infty} |\varphi_{\pm}(t,h)|^k=\dfrac{1}{2}\cdot \dfrac{|\varphi_{\pm}(t,h)|^2}{1-|\varphi_{\pm}(t,h)|}\leqslant |\varphi_{\pm}(t,h)|^2.
\end{eqnarray*}
Thus we get the estimate
\begin{eqnarray*}
	|\delta(t,h)|\leqslant \tfrac{1}{2}|\varphi_{\pm}(t,h)|^2+\tfrac{1}{2} |\rho(t,h)|,\quad t\geqslant0,\,\, h>0,
\end{eqnarray*}
under the condition $\bigl| \varphi_{\pm}(t,h)\bigr|\leqslant \tfrac{1}{2}$.

Let  us  choose the sequence $(h_l)_{l\in\N}$ as in the statement of the theorem, i.e.  $h_l\to 0$ as $l\to\infty$, \eqref{th_suff_phi_b} and \eqref{th_suff_phi_phi2} hold. We also arbitrarily fix $t>0$ and define the integers $n_{t,l}\defeq \lfloor t/ h_l\rfloor\in \N_0$, $l\in\N$.  So we have $\varphi_{\pm}(kh_l,h_l)\to 0$, $l\to\infty$, and, in particular,  $| \varphi_{\pm}\bigl(kh_l,h_l)\bigr|\leqslant \tfrac{1}{2}$ for all sufficiently large $l$ uniformly over $k=0,\ldots, n_{t,l}$. This is seen due to the expression $\varphi_{\pm}$ through the values of  $f$ on account of the uniform continuity of  $f$ and of the assumption $f(t) \ne 0$, $t\in\R$.  So for all sufficiently large $l \in\N$ we have
\begin{eqnarray*}
	\bigl|\delta(kh_l,h_l)|\leqslant \tfrac{1}{2}|\varphi_{\pm}(kh_l, h_l )|^2+ \tfrac{1}{2} |\rho(kh_l,h_l)\bigr|\quad \text{for all}\quad k=0,\ldots, n_{t,l}.
\end{eqnarray*}
This estimate together with  \eqref{th_suff_phi_b} and \eqref{th_suff_phi_phi2}  yield \eqref{th_crit_psi_cond_a}. Then, by Theorem \ref{th_crit_psi}, all assertions of Theorem \ref{th_suff_phi} hold. \quad $\Box$\\

\textbf{Proof of Theorem \ref{th_SecDerivLnf}.} a) Under the assumptions of the theorem, we first show that for any $t\in \R$ there exists $(\Ln f)'(t)$. Let us fix $t\in\R$ and consider the relation
\begin{eqnarray}\label{def_diff_Lnfth}
	\dfrac{\Ln f(t+h)-\Ln f(t)}{h}=i \gamma +\int_{\R}\biggl( \dfrac{e^{i(t+h)x}-e^{itx}}{h}-i \sin(x)\biggr)\tfrac{1+x^2}{x^2}\dd G(x),
\end{eqnarray}
where $0<|h|<1$. We will propose majorant for  the integrated function.   Let $r$ be a positive real number. For any $x\in\R$ satisfying  $|x|> r$ we have the estimate
\begin{eqnarray*}
	 \biggl|\dfrac{e^{i(t+h)x}-e^{itx}}{h}-i \sin(x)\biggr|\leqslant  \biggl|\dfrac{e^{i(t+h)x}-e^{itx}}{h}\biggr| + |i\sin(x)| = \biggl|\dfrac{e^{ihx}-1}{h}\biggr| + |\sin(x)|.
\end{eqnarray*}
In the last sum, each term is bounded by $|x|$. So we have 
\begin{eqnarray*}
	\biggl|\dfrac{e^{i(t+h)x}-e^{itx}}{h}-i \sin(x)\biggr|\leqslant  2|x|, \quad |x|> r. 
\end{eqnarray*}
For the case $|x|\leqslant r$ we will use the representation
\begin{eqnarray*}
	\dfrac{e^{i(t+h)x}-e^{itx}}{h}-i \sin(x)=\dfrac{e^{i(t+h)x}-e^{itx}-ihx}{h}-i (\sin(x)-x). 
\end{eqnarray*}
So we have the inequality
\begin{eqnarray*}
	\biggl|\dfrac{e^{i(t+h)x}-e^{itx}}{h}-i \sin(x)\biggr|\leqslant \biggl|\dfrac{e^{i(t+h)x}-e^{itx}-ihx}{h}\biggl| + | \sin(x)-x |.
\end{eqnarray*}
For the last term we have
\begin{eqnarray*}
	| \sin(x)-x |\leqslant \dfrac{|x|^3}{3!}\leqslant \dfrac{r x^2}{3!},  \quad |x|\leqslant r.
\end{eqnarray*}
Next, we represent
\begin{eqnarray*}
	\dfrac{e^{i(t+h)x}-e^{itx}-ihx}{h} =   \dfrac{ix}{h}\int_{t}^{t+h} (e^{iux}-1)\dd u =  \dfrac{ix}{h}\int_{0}^{h} (e^{i(v+t)x}-1)\dd v.
\end{eqnarray*}
Since $\bigl| e^{i(v+t)x}-1 \bigr|\leqslant \bigl|  (v+t)x\bigr|\leqslant \bigl(|v|+|t|\bigr)| x |$, we get
\begin{eqnarray*}
\Biggl| \int_{0}^{h} (e^{i(v+t)x}-1)\dd v\Biggr| \leqslant \Biggl| \int_{0}^{h}\bigl| e^{i(v+t)x}-1 \bigr|\dd v\Biggr| \leqslant \Biggl| \int_{0}^{h} \bigl(|v|+|t|\bigr)| x |\dd v\Biggr| = | x |\int_{0}^{|h|} \bigl(v+|t|\bigr)\dd v.
\end{eqnarray*}
Therefore
\begin{eqnarray*}
	 \biggl|\dfrac{e^{i(t+h)x}-e^{itx}-ihx}{h}\biggr| = \Biggl|  \dfrac{x}{h}\int_{0}^{h} \bigl( e^{i(v+t)x}-1 \bigr) \dd u \Biggr|\leqslant  \dfrac{x^2}{|h|}\int_{0}^{|h|}\bigl(v+|t|\bigr)  \dd v  \leqslant  \bigl( |h|+|t|   \bigr) x^2.
\end{eqnarray*}
Since $0<|h|<1$, we have
\begin{eqnarray*}
	\biggl|\dfrac{e^{i(t+h)x}-e^{itx}}{h}-i \sin(x)\biggr|\leqslant C_{t,r} x^2, \quad |x|\leqslant r,
\end{eqnarray*}
with $C_{t,r}= 1+|t|+\tfrac{r}{3!}$. 

The obtained upper estimates are integrable with $\tfrac{1+x^2}{x^2}$ by $|G|(\cdot)$ over the coresponding sets :
\begin{eqnarray*}
	\int_{|x|> r} 2|x|\,\tfrac{1+x^2}{x^2}\dd |G|(x)=2\int_{|x|> r} \tfrac{1+x^2}{|x|}\dd |G|(x)\leqslant \dfrac{2}{r}\int_{|x|> r} (1+x^2)\dd |G|(x)   <\infty,
\end{eqnarray*}
and
\begin{eqnarray*}
	\int_{|x|\leqslant r} C_{t,r} x^2\,\tfrac{1+x^2}{x^2}\dd |G|(x)=C_{t,r}\int_{|x|\leqslant r} (1+x^2)\dd |G|(x)<\infty.
\end{eqnarray*}
Thus the integrand in \eqref{def_diff_Lnfth} has the integrable majorant. 

Next, we observe that the integrand in \eqref{def_diff_Lnfth}  converges for any  $x\ne 0$:
\begin{eqnarray}\label{conc_lim_lnf'}
	\lim_{h\to 0}\biggl( \dfrac{e^{i(t+h)x}-e^{itx}}{h}-i \sin(x)\biggr)\tfrac{1+x^2}{x^2}=\bigl(ix e^{itx}-i \sin(x)\bigr)\tfrac{1+x^2}{x^2}.
\end{eqnarray}
Defining the right-hand side by continuity at $x=0$, it is seen that the equality is true for this point too. Indeed, due to the well-known expansions
\begin{eqnarray*}
	e^{itx}=1+itx+ O(x^2), \quad  \sin (x)=x-O(x^3),\quad x\to 0,
\end{eqnarray*}
we have
\begin{eqnarray*}
	\bigl(ix e^{itx}-i \sin(x)\bigr)\tfrac{1+x^2}{x^2}= \bigl((ix)^2 t+O(x^3)\bigr)\tfrac{1+x^2}{x^2}\to -t,\quad x\to 0.
\end{eqnarray*}
By the following well-known convention
\begin{eqnarray*}
	\Bigl(\bigl( e^{itx}-1-it \sin( x)\bigr)\tfrac{1+x^2}{x^2}\Bigr)\biggl|_{x=0}\defeq\lim\limits_{x\to 0}\Bigl(\bigl( e^{itx}-1-it \sin( x)\bigr)\tfrac{1+x^2}{x^2}\Bigr)=-\dfrac{t^2}{2},
\end{eqnarray*}
we have
\begin{eqnarray*}
	\lim_{h\to 0}\biggl( \dfrac{e^{i(t+h)x}-e^{itx}}{h}-i \sin(x)\biggr)\tfrac{1+x^2}{x^2}= \lim_{h\to 0}\biggl\{ \dfrac{1}{h}\biggl(-\dfrac{(t+h)^2}{2}+ \dfrac{t^2}{2} \biggr) \biggr\} = \lim_{h\to 0}\biggl\{ \dfrac{1}{h}\biggl(-th - \dfrac{h^2}{2}\biggr) \biggr\}=-t.
\end{eqnarray*}
Thus \eqref{conc_lim_lnf'} holds for all $x\in\R$.

Due  to the Jordan decomposition, we decompose $G(x)=G^+(x)-G^-(x)$, $x\in\R$,  with some non-decreasing functions $G^+$ and $G^-$ from $\BV$ that satisfy
\begin{eqnarray*}
	\int_{\R} (1+x^2)\dd |G^+|(x)<\infty,\quad \text{and}\quad  \int_{\R} (1+x^2)\dd |G^-|(x)<\infty.
\end{eqnarray*}
For example, $G^+(x)\defeq \tfrac{1}{2} \bigl(|G|(x)+G(x)\bigr)$ and $G^-(x)\defeq \tfrac{1}{2} \bigl(|G|(x)-G(x)\bigr)$, $x\in\R$. So we write
\begin{eqnarray*}
	\int_{\R}\biggl( \dfrac{e^{i(t+h)x}-e^{itx}}{h}-i \sin(x)\biggr)\tfrac{1+x^2}{x^2}\dd G(x)&=& \int_{\R}\biggl( \dfrac{e^{i(t+h)x}-e^{itx}}{h}-i \sin(x)\biggr)\tfrac{1+x^2}{x^2}\dd G^+(x)\\
	&&{}\qquad-\int_{\R}\biggl( \dfrac{e^{i(t+h)x}-e^{itx}}{h}-i \sin(x)\biggr)\tfrac{1+x^2}{x^2}\dd G^-(x).
\end{eqnarray*}
Due to \eqref{conc_lim_lnf'} and the proved domination, we now apply Lebesgue's dominated convergence theorem to the integrals in the right-hand side as $h\to 0$:
\begin{eqnarray*}
	\lim_{h\to 0}\int_{\R}\biggl( \dfrac{e^{i(t+h)x}-e^{itx}}{h}-i \sin(x)\biggr)\tfrac{1+x^2}{x^2}\dd G(x)
	&=&\int_{\R}\bigl(ix e^{itx}-i \sin(x)\bigr)\tfrac{1+x^2}{x^2}\dd G^+(x)\\
	&&{}\qquad-\int_{\R}\bigl(ix e^{itx}-i \sin(x)\bigr)\tfrac{1+x^2}{x^2}\dd G^-(x)\\
	&=&\int_{\R}\bigl(ix e^{itx}-i \sin(x)\bigr)\tfrac{1+x^2}{x^2}\dd G(x).	
\end{eqnarray*}
Thus, according to formula \eqref{def_diff_Lnfth}, for any $t\in\R$ we have 
\begin{eqnarray*}
	(\Ln f)'(t)=\lim_{h\to 0}\dfrac{\Ln f(t+h)-\Ln f(t)}{h}
	=i \gamma +\int_{\R}\bigl(ix e^{itx}-i \sin(x)\bigr)\tfrac{1+x^2}{x^2}\dd G(x).
\end{eqnarray*}
We now prove that $(\Ln f)''(t)$ exists for any $t\in\R$ and that it is represented by formula \eqref{th_SecDerivLnf_2deriv}. Let us fix $t\in\R$ and consider the relation
\begin{eqnarray*}
	\dfrac{(\Ln f)'(t+h)-(\Ln f)'(t)}{h}=\int_{\R}\biggl( ix \cdot\dfrac{e^{i(t+h)x}-e^{itx}}{h}\biggr)\tfrac{1+x^2}{x^2}\dd G(x),\quad  h\in\R\setminus\{0\}.
\end{eqnarray*}
Observe that
\begin{eqnarray*}
	\lim_{h\to 0}\Biggl[\biggl( ix \cdot\dfrac{e^{i(t+h)x}-e^{itx}}{h}\biggr)\dfrac{1+x^2}{x^2}\Biggr]= ix e^{itx}  \lim_{h\to 0}\biggl(  \dfrac{e^{ihx}-1}{h}\biggr)\cdot\dfrac{1+x^2}{x^2} =(ix)^2 e^{itx}\,  \dfrac{1+x^2}{x^2}=-e^{itx} (1+x^2),
\end{eqnarray*}
with the following domination uniformly over $h\in\R\setminus\{0\}$:
\begin{eqnarray*}
	\Biggl| \biggl( ix \cdot\dfrac{e^{i(t+h)x}-e^{itx}}{h}\biggr)\dfrac{1+x^2}{x^2} \Biggr|= |x| \cdot\biggl|\dfrac{e^{ihx}-1}{h}\biggr|\cdot \dfrac{1+x^2}{x^2} \leqslant |x| \cdot  |x|\cdot \dfrac{1+x^2}{x^2}=1+x^2.
\end{eqnarray*}
Due to \eqref{th_SecDerivLnf_int}, the analogous application of Lebesgue's dominated convergence theorem  is valid and  it yields 
\begin{eqnarray*}
	\lim_{h\to 0}\int_{\R}\biggl( ix \cdot\dfrac{e^{i(t+h)x}-e^{itx}}{h}\biggr)\tfrac{1+x^2}{x^2}\dd G(x)=-\int_{\R}e^{itx}(1+x^2)\dd G(x).
\end{eqnarray*}
Thus for any $t\in\R$ we obtain the required equality
\begin{eqnarray*}
	(\Ln f)''(t)=	\lim_{h\to 0} \dfrac{(\Ln f)'(t+h)-(\Ln f)'(t)}{h} = -\int_{\R}e^{itx}(1+x^2)\dd G(x).
\end{eqnarray*}

b)  We now assume that \eqref{th_SecDerivLnf_2deriv}  holds with some function $G\in\BV$ that satisfies  \eqref{th_SecDerivLnf_int}. We first will show that the function $(\Ln f)''$ is uniformly continuous on $\R$. For any  $h \in \R$ we have
\begin{eqnarray*}
\sup_{t\in\R}\bigl|(\Ln f)''(t+h) -(\Ln f)''(t)\bigr|&=&\sup_{t\in\R} \biggl| \int_{\R} (e^{i(t+h)x}-e^{itx})(1+x^2)\dd G(x)\biggr| \\
&\leqslant&\sup_{t\in\R} \int_{\R} \bigr|e^{i(t+h)x}-e^{itx}\bigl|(1+x^2)\dd |G|(x)\\
&\leqslant& \int_{\R} |e^{ihx}-1|(1+x^2)\dd |G|(x).
\end{eqnarray*}
Hence for any $r>0$ the following estimate holds
\begin{eqnarray*}
	 \sup_{t\in\R}\bigl|(\Ln f)''(t+h) -(\Ln f)''(t)\bigr|&\leqslant&\int_{|x|\leqslant r} |e^{ihx}-1|(1+x^2)\dd |G|(x)+ \int_{|x|>r} |e^{ihx}-1|(1+x^2)\dd |G|(x)\\
	 &\leqslant&\int_{|x|\leqslant r} |hx|(1+x^2)\dd |G|(x)+ \int_{|x|>r}\bigl(|e^{ihx}|+ 1  \bigr) (1+x^2)\dd |G|(x)\\
	 &\leqslant& |h| r\int_{|x|\leqslant r} (1+x^2)\dd |G|(x)+ 2\int_{|x|>r} (1+x^2)\dd |G|(x).
\end{eqnarray*}
Let us fix arbitrary $\e>0$. Due to \eqref{th_SecDerivLnf_int}, we first chose $r=r_\e$ such that
\begin{eqnarray*}
	\int_{|x|>r_\e} (1+x^2)\dd |G|(x)\leqslant \e,
\end{eqnarray*}
and after this we choose $h_\e>0$ such that
\begin{eqnarray*}
	h_\e r_\e \int_{|x|\leqslant r_\e} (1+x^2)\dd |G|(x)\leqslant \e.
\end{eqnarray*}
Then for any $h\in\R$ satisfying $|h|\leqslant h_\e$ we obtain
\begin{eqnarray*}
	\sup_{t\in\R}\bigl|(\Ln f)''(t+h) -(\Ln f)''(t)\bigr| \leqslant h_\e r_\e\int_{|x|\leqslant r} (1+x^2)\dd |G|(x)+ 2\int_{|x|>r_\e} (1+x^2)\dd |G|(x)\leqslant 3\e.
\end{eqnarray*}
This means that $(\Ln f)''$ is uniformly continuous on $\R$.

Due to the proved continuity of $(\Ln f)''$, we can write
\begin{eqnarray*}
	(\Ln f)'(s) =(\Ln f)'(0)+\int_{0}^s  (\Ln f)''(u)\dd u,\quad s\in \R.
\end{eqnarray*}
Since $\Ln f(0)=0$ and $(\Ln f)'$ must be continuous, we have
\begin{eqnarray*}
	\Ln f(t) = \int_{0}^t  (\Ln f)'(s)\dd s,\quad t\in \R.
\end{eqnarray*}
Therefore
\begin{eqnarray*}
	\Ln f(t) = \int_{0}^t \biggl( (\Ln f)'(0)+\int_{0}^s  (\Ln f)''(u)\dd u \biggr) \dd s=   t \cdot (\Ln f)'(0)+\int_{0}^t \biggl(\int_{0}^s  (\Ln f)''(u)\dd u \biggr) \dd s,\quad t\in \R.
\end{eqnarray*}
Let us fix $t \in \R$ and $h>0$. We write
\begin{eqnarray*}
	\psi(t,h) &=&t \cdot (\Ln f)'(0)+\int_{0}^{t} \biggl(\int_{0}^s  (\Ln f)''(u)\dd u \biggr) \dd s\\
	&&{}-\dfrac{1}{2}\,(t-h) \cdot (\Ln f)'(0)-\dfrac{1}{2}\int_{0}^{t-h} \biggl(\int_{0}^s  (\Ln f)''(u)\dd u \biggr) \dd s\\
	&&{}-\dfrac{1}{2}\,(t+h) \cdot (\Ln f)'(0)-\dfrac{1}{2}\int_{0}^{t+h} \biggl(\int_{0}^s  (\Ln f)''(u)\dd u\biggr) \dd s.
\end{eqnarray*}
Thus we get
\begin{eqnarray*}
\psi(t,h)=-\dfrac{1}{2}\int_{t}^{t+h} \biggl(\int_{0}^s  (\Ln f)''(u)\dd u \biggr) \dd s -\dfrac{1}{2} \int_{t}^{t-h} \biggl(\int_{0}^s  (\Ln f)''(u)\dd u \biggr) \dd s.
\end{eqnarray*}
On account of \eqref{th_SecDerivLnf_2deriv}, we  have
\begin{eqnarray*}
	\int_{t}^{t\pm h} \biggl(\int_{0}^s  (\Ln f)''(u)\dd u \biggr) \dd s &=& -\int_{t}^{t\pm h} \biggl(\int_{0}^s \biggl(  \int_{\R} e^{iux}(1+x^2)\dd G(x)\biggr) \dd u \biggr) \dd s \\
	& =& -\int_{\R}\biggl(\int_{t}^{t\pm h} \biggl(\int_{0}^s    e^{iux}\dd u \biggr) \dd s\biggr) (1+x^2)\dd G(x).
\end{eqnarray*}
It is easily seen that the changing of the order of integration is valid  due to the boundedness of  $e^{iux}$ for $u,x \in \R$.  So we get
\begin{eqnarray*}
	\psi(t,h)=\dfrac{1}{2}\int_{\R}\Biggl[\int_{t}^{t+ h} \biggl(\int_{0}^s    e^{iux}\dd u \biggr) \dd s+ \int_{t}^{t- h} \biggl(\int_{0}^s    e^{iux}\dd u \biggr) \dd s\Biggr] (1+x^2)\dd G(x).
\end{eqnarray*}
Observe that
\begin{eqnarray*}
	\int_{t}^{t\pm h} \biggl(\int_{0}^s    e^{iux}\dd u \biggr) \dd s= \int_{t}^{t\pm h} \biggl(\dfrac{    e^{isx}-1 }{ix}\biggr) \dd s=\dfrac{    e^{i(t\pm h)x}-e^{itx}}{(ix)^2} -\dfrac{  \pm h }{ix}.
\end{eqnarray*}
Hence
\begin{eqnarray*}
	\psi(t,h)
	&=& \dfrac{1}{2}\int_{\R}\biggl(	\dfrac{    e^{i(t+ h)x}-e^{itx}}{(ix)^2}+ 	\dfrac{  e^{i(t-h)x}-e^{itx}}{(ix)^2}\biggr) (1+x^2)\dd G(x)\\
	&=& \int_{\R}e^{itx}\cdot	\dfrac{    e^{ihx}+e^{-ihx}-2}{2(ix)^2}\, (1+x^2)\dd G(x)\\
	&=&\int_{\R}e^{itx}\cdot	\dfrac{  \cos (hx)-1}{(ix)^2}\,	(1+x^2)\dd G(x).
\end{eqnarray*}
Thus for any $t \in \R$ and $h>0$ we obtain
\begin{eqnarray*}
	\psi(t,h)=\int_{\R}e^{itx}\bigl(1-\cos(hx)\bigr)\tfrac{1+x^2}{x^2}\dd G(x).
\end{eqnarray*}
Using Corollary  \ref{co_th_crit_psi}, we come to the required assertion.

If we additionally suppose that $f$ is the characteristic function of some distribution function $F$, then $F$ is quasi-infinitely divisible and, consequently, it is rational-infinitely divisible, i.e. $F\in\Q$ (see Introduction). \quad $\Box$\\

\section{Acknowledgments}

This research was supported by the Ministry of Science and Higher Education of the Russian Federation, agreement 075-15-2022-289 date 06/04/2022.

\end{document}